\documentclass[10pt]{article}
\usepackage{amsmath, amssymb}
\usepackage{latexsym}
\usepackage{mathtools}
\usepackage{graphicx}
\usepackage{color}
\usepackage{url}
\numberwithin{equation}{section}
\DeclareMathAlphabet{\itbf}{OML}{cmm}{b}{it}

\newcommand{\RR}{\mathbb{R}}

\newcommand{\NN}{\mathbb{N}}

\newcommand{\ds}{\displaystyle}
\newcommand{\no}{\nonumber}
\newcommand{\ri}{\rightarrow}

\newcommand{\q}{\quad}

\newcommand{\bx}{{\itbf x}}
\newcommand{\bv}{{\itbf v}}
\newcommand{\bg}{{\itbf g}}
\newcommand{\bh}{{\itbf h}}
\newcommand{\bev}{{\itbf e}}

\newcommand{\bu}{{\itbf u}}
\newcommand{\bn}{{\itbf n}}
\newcommand{\by}{{\itbf y}}

\newcommand{\bi}{\begin{itemize}}
\newcommand{\ei}{\end{itemize}}

\newcommand{\bH}{{\itbf H}}

\newcommand{\be}{\begin{eqnarray}}
\newcommand{\ee}{\end{eqnarray}}

\newcommand{\ben}{\begin{eqnarray*}}
\newcommand{\een}{\end{eqnarray*}}
\def\ds{\displaystyle}

\newcommand\ov{\overline}

\newtheorem{lem}{Lemma}[section]
\newtheorem{prop}{Proposition}[section]

\newtheorem{thm}{Theorem}[section]

\newtheorem{definition}{Definition}[section]

\newcommand{\bea}{\begin{eqnarray*}}
\newcommand{\eea}{\end{eqnarray*}}
\newcommand{\bean}{\begin{eqnarray}}
\newcommand{\eean}{\end{eqnarray}}
\newcommand{\p}{\partial}
\newcommand{\f}{\frac}
\newcommand{\s}{\sqrt}
\newcommand{\di}{\mbox{div }}
\newcommand{\aaa}{\mbox{$[$}}
\newcommand{\bbb}{\mbox{$]$}}

\begin{document}
% corrected typos post submission
%

\title{A stochastic approach to reconstruction of faults in elastic half space}

\author{ Darko
Volkov, \thanks{\footnotesize D. Volkov is supported  by
a Simons Foundation Collaboration Grant.} \thanks{Department of Mathematical Sciences,
Worcester Polytechnic Institute, Worcester, MA 01609.
} \and Joan Calafell Sandiumenge \thanks{Heat and Mass Technological Center (CTTC), Technical University
of Catalonia (UPC), Colom 11, 08222 Terrassa (Barcelona), Spain.}}
%joancs@cttc.upc.edu
%Heat and Mass Technological Center (CTTC). Technical University
%of Catalonia (UPC), Colom 11, 08222 Terrassa (Barcelona), Spain

\maketitle

\begin{abstract}
We introduce in this study an algorithm 
for the  imaging of faults and of slip fields on those faults.
 %to simultaneously reconstruct the geometry of a fault 
%and the slip field on that fault based on surface measurements of displacement fields.
The physics of this problem are modeled  using the equations of linear elasticity.
%and the data for the fault inverse problem consists of surface measurements of 
%displacements. 
We define a regularized functional to be minimized for building the image.
We first prove that the minimum of that functional
%solution to this minimization problem 
converges to the   unique solution
of the related fault inverse problem.
Due to inherent uncertainties in measurements,
rather than seeking a deterministic solution to the fault inverse problem,
 we then consider a Bayesian approach. 
In this approach the geometry of the fault is assumed to be planar, it can thus be modeled 
by a three dimensional random variable whose probability density has to be determined knowing
surface measurements. 
The randomness involved in the unknown slip is teased out by assuming independence of the priors,
and we show how
 the regularized error functional introduced earlier can be used to recover the probability density of the geometry parameter. 
%The standard deviation of the slip (or any moment of the slip) is computed
%in a separate step where the geometry of the fault is fixed (for example, the most likely geometry or the mean geometry). 
% to solving the fault inverse problem.
The advantage of the Bayesian  approach is that we obtain a way of quantifying uncertainties 
as part of  our final answer. 
On the downside, this 
 approach leads to a very large computation since the slip is unknown.
To contend with the size of this computation
we  developed  an algorithm for the numerical solution to the stochastic minimization problem which can be easily
implemented on  a parallel multi-core  
	platform and we discuss techniques aimed at saving on computational time. 
%Although it involves solving a large number of PDEs, we explain how our algorithm 
	%These techniques rely on saving some pre-computed terms which are common to all the processes
	%and applying adequate numerical linear algebra methods.
After showing how this algorithm performs on simulated data,
we apply it to  measured data. The data was recorded during a slow slip event 
in  Guerrero, Mexico.
%We study in this paper a half space linear elasticity model for
%surface displacements caused by slip along underground faults. 
%This model is commonly used in geophysics. 
%We prove uniqueness
%of the fault location and (piecewise planar) geometry and of the slip field for a given
%surface displacement field.
%We then introduce a reconstruction algorithm for the realistic case where only 
%a finite number of surface measurements are available.
%
%Since this is a well studied subduction zone, it is possible to compare
%our inferred fault geometry to other reconstructions (obtained using different techniques),
%found in the literature.   
\end{abstract}

\bigskip

%\noindent {\footnotesize Mathematics Subject Classification
%(MSC2000): 35R30, 35B30}

%\noindent {\footnotesize Keywords: eddy current imaging, induction data, classification, recognition,
%invariant shape descriptors}

\section*{Acknowledgements}
Results in this paper were obtained in part using a high-performance computing system acquired through NSF MRI grant DMS-1337943 to WPI.\\
D. Volkov is supported  by
a Simons Foundation Collaboration Grant.

\section{Introduction}
%Recently published Nature paper \\
%
Subduction zones around the world 
are periodically prone to devastating earthquakes.
%ose a major threat in terms of seismic hazard. 
%This is where
%the largest earthquakes occur, sometimes followed by devastating tsunamis. 
The 2011 Tohoku
Oki earthquake in Japan was a stark reminder of that occurrence.
Experts are now warning the public and policy makers that in North America,
 the Pacific Northwest is in great danger of being struck by a massive earthquake
in the near future partly because there is strong evidence
of such events in the recent past \cite{atwater1995summary}. Of course it is not possible at this stage to say when  this may happen again.
A better knowledge of the structure of the subduction zone 
and the interface between oceanic and continental crusts
in the Pacific Northwest, together with an assessment of the mechanical stress budget
will undoubtedly help geophysicists make progress in predictive skills.
Deformations in the vicinity of  major subduction zones 
have been continuously recorded for some time using
%is scrutinized with 
 geodetic networks (GPS, tiltmeters) as well as broadband  seismological
networks. 
%These sophisticated instrumentation 
%networks have led to  uncovering  previously unknown phenomena such as 
% Non Volcanic Tremors (NVT) and Low
%Frequency Earthquakes (LFE) \cite{Obara2002}, which were first discovered in Japan, and are now recognized
%everywhere along subduction zones and also on continental faults 
%\cite{GuilhemNadeau2012, Guilhem2010}.
 GPS networks 
%(now complemented by InSAR observations)
 have revealed the
existence of periods of reversed motion relative to the interseismic motions in many subduction
zones worldwide \cite{Dragert2004,Dragert2001}. 
These reversals of  movement are
interpreted as aseismic slow slip events (SSE) occurring 
 deep beneath the subduction zone below the
locked seismogenic zone. 
%These SSE show an important variability in terms of duration (from days
%to years), recorded surface displacements (a few millimeters to a few centimeters), and recurrence
%time (see \cite{SchwartzRokosky} for a review). \\
%These SSE's have changed the current understanding of  seismic cycles of  subduction zones
%\cite{araki2017recurring}.
%Depending on their precise location, extent and amount of slip, they may have opposite
%implications on  stress budgets of  subduction zones, and therefore also on the occurrence
% of large  earthquakes near subduction zones. 
%A detailed knowledge of the subduction geometry
% and
%the slip dynamics  is therefore of prime interest, but has proved difficult to obtain.
Prior to two recent  studies \cite{volkov2017determining, volkov2017reconstruction}, 
the geometry profiles of  subduction zones have been derived for the most part
  from seismology 
and are therefore poorly constrained. In a separate step, using these geometry profiles, 
investigators have studied
slip distributions for SSE's from GPS time series, occasionally augmented by InSAR data.\\
In this paper, we introduce and analyze error functionals  for the   reconstruction of fault 
geometries
based on surface measurements of displacement fields, and we derive a stochastic inversion procedure which relies 
on these functionals.
The physics of our problem are modeled  using the equations of linear elasticity
and the data for the fault inverse problem consists of measurements of surface   
displacements. 
%which is  common in the geophysics community. 
Evidently, GPS surface measurements are inherently tainted by errors. There are also
errors due to using a PDE model, which of course can only give a simplified sketch
of complex geophysical processes occurring in subduction zones.
In this paper we take into account these errors 
by seeking to determine probability densities for the 
geometry 
thanks to a Bayesian formulation. We assume that the fault is planar which is a common assumption in
geophysics, at least for the active part of a fault (the part where the slip is non zero)  during a slow slip event. With this assumption only the joint probability density of three scalar
parameters giving the equation of the plane containing the fault has to be determined, thanks to the assumption that the geometry parameters and the slip field on the fault are independent.\\
% through the use
%of the Bayesian formula.\\
This paper is organized as follows. In section \ref{Mathematical}
we introduce a mathematical formulation for modeling slow slip events on faults.
This model relies on the equations of linear isotropic elasticity in  half space. We review existence and uniqueness results
for the forward problem and a uniqueness result for the fault inverse problem.
This result ensures that if surface displacement fields are known on an open set of the 
top boundary then it is possible to reconstruct the fault and the slip on that fault.
In section \ref{Continuous regularized reconstruction} we introduce a regularized functional for 
reconstructing faults and slip fields
from surface measurements. We prove that as the regularization constant $C$ tends to zero,
the reconstructed profile and slip field obtained by minimizing that functional 
converge to the actual profile and slip that produced the surface displacements.
Let us point out here that this result is not trivial since, although this inverse problem
is linear in the slip field, it is {\sl non -linear} in the geometry of the fault.
In section \ref{finite recon} we define another reconstruction functional which 
involves only  a finite number of surface measurements and slip fields in a finite dimensional
subspace. 
In that  case it is not possible to invoke the uniqueness result for the inverse problem proved 
in earlier work. However, we are able to prove that if the number of measurement points is 
sufficiently large,
if the subspace over which this functional is minimized is large enough, 
 and if the regularization constant is sufficiently small, then 
the solution to this discrete minimization  problem  can be  
 arbitrarily close to the actual profile and slip that produced the surface displacements.\\
%minimum of the function $f^{disc}$ occurs 
%arbitrarily close to $(\tilde{a},\tilde{b},\tilde{d})$.
In section \ref{Stochastic modeling} %247
we take into consideration that
 the number $N$ of surface displacement measurements  is low
and these measurements are uncertain, 
so rather than seeking a deterministic solution to the fault inverse problem,
 we  consider a Bayesian approach to solving the fault inverse problem.
%Our formulation is closely related to the one considered  by Kaipio and Somersalo 
%in section 3.2.2 of \cite{kaipio2006statistical}. 
As customary  in Bayesian modeling of inverse problems,
 the difference between 
measured data 
and predicted surface displacements for a given geometry of the fault  and a given slip
is assumed to be a
 Gaussian random variable
with mean zero.
% and  given covariance which can be estimated from the data.
%, multiplied by a Gaussian 
%that penalizes highly oscillatory slip fields.
The random slip field is also assumed to be  Gaussian and we make the assumption that
the prior geometry parameters and the prior slip field are independent. 
Thanks to this independence assumption, recovering only the probability density of the geometry
parameters becomes a computationally tractable problem, albeit by use of advanced computational techniques discussed further.
%In a second step, the probability density of the geometry knowing the measurements
%is computed thanks to the Bayesian formula. 
As  further motivation for this Bayesian approach, 
we  prove in section 5 that the recovered probability density of the geometry parameters
tends to zero for all geometry parameters different from those of the true profile as the number
of measurement points grows large and the variance of the measurements and the regularization 
parameter for the reconstructed slip field become small.\\
% (and this convergence can be made uniform on adequately defined sets). \\
%After teasing out the randomness of the slip field, our  formulation still leads to a very large computation since each candidate for 
%a possible geometry spans its own set of calculations.
%We thus developed an algorithm that takes into account that some terms remain unchanged 
%between geometries and can therefore be pre-computed. 
%   Instead of solving a large linear system for each geometry we explain how we take advantage
%	of Woodbury's formula: it allows us to solve  much smaller linear systems for 
	%each geometry and to use a pre-computed inverse of a large matrix.
	Our proposed algorithm is amenable  to implementation on a parallel multi-core  computational
	platform. The combination of relevant linear algebra techniques
	 and parallel implementation  led to great savings in computational time.
	In section \ref{Numerical results} %268
we apply this reconstruction algorithm to the case of the 2007 Guerrero, Mexico SSE.
 We first examine three test cases with numerically generated data for the inverse problem.
The surface points for the test cases are the same as those where geophysicists 
sampled real world measurements. All length scales and noise levels have same order
of magnitude 
 as those
observed in the real world. Different geometries are considered and in one case we add 
a systematic error due to imperfections in the model.
Our last numerical computation involves real world measurements and results
 in the reconstruction 
of the
part of the subduction interface beneath the  Guerrero region which was active 
during the 2007 SSE. 
In this last simulation the only benchmarks for our calculation are 
geometries estimated by other authors (in most cases, based on other physical processes).
We observe that many of the profiles found by other authors fall in the plus or minus one 
standard deviation envelope of  the profile derived in this present study.
%Our reconstructed ge 

\section{Mathematical model and uniqueness result} \label{Mathematical}
%We denote the natural basis of $\RR^{3}$
%$\bev_1, \bev_2, \bev_3$ throughout  this paper.

\subsection{Forward problem}
Using the standard rectangular coordinates $\bx= (x_1, x_2, x_3)$ of $\RR^{3}$,
%A point $(x,y,z)$ will denote a point 
we define $\RR^{3-}$ to be the open half space $x_3<0$.
%Let  $\Gamma$ be a bounded open surface
%with smooth boundary included in the plane $\{ x_3=0 \}$.
The derivative in the $i$-th coordinate will be denoted by $\p_i$.
In this paper 
we only consider the case of linear, homogeneous, isotropic elasticity; 
the two Lam\'e constants $\lambda$ and $\mu$ will be two positive constants.
For a vector field $\bu = (u_1, u_2, u_3)$, 
the stress  and strain tensors will be denoted as follows,
\bea
\sigma_{ij}(\bu) = \lambda \, \di \bu \, \delta_{ij} + \mu \, (\p_i u_j + \p_j u_i ), \\
\epsilon_{ij}(\bu) = \f12 (\p_i u_j + \p_j u_i ),
\eea
and the  stress vector in the normal 
direction $\bn$ will be denoted by
\bea
T_n \bu = \sigma (\bu) \bn.
\eea

Let $\Gamma$ be a Lipschitz open surface which is strictly included in $\RR^{3-}$.
Let $\bu$ be the displacement field solving
\bean
\mu \Delta \bu+ (\lambda+\mu) \nabla \di \bu= 0  \mbox{ in } \RR^{3-} 
\setminus \Gamma \label{uj1}, \\
 T_{\bev_3} \bu =0 \mbox{ on the surface } x_3=0 \label{uj2}, \\
 T_{\bn} \bu  \mbox{ is continuous across } \Gamma \label{uj3}, \\
 \aaa  \bu\bbb =\bg \mbox{ is a given jump across } \Gamma , \label{uj3andhalf}\\
\bu (\bx) = O(\f{1}{|\bx|^2}),  \nabla \bu (\bx) = O(\f{1}{|\bx|^3}), \mbox{ uniformly as } 
|\bx| \rightarrow \infty,
\label{uj4}
\eean
where %$\bx=(x_1, x_2, x_3)$ is any point in $\RR^{3-}$ and
$\bev_3$ is the vector $(0,0,1)$.
In \cite{volkov2017reconstruction}, we defined the functional space 
$\cal{V}$ of vector fields $\bu $ defined in $\RR^{3-}\setminus{\ov{\Gamma}}$
such that $\nabla \bu $ and $\ds \f{\bu}{(1+ r^2)^{\f12}}$ are in $L^2(\RR^{3-}\setminus{\ov{\Gamma}})$.
 Let $\p D$ be a closed Lipschitz surface containing $\Gamma$. We define the Sobolev space 
  $\tilde{H}^{\f12}(\Gamma)^2$ 
 to be the set of restrictions to $\Gamma$ of tangential fields  in
   $H^{\f12}(\p D)^2$ supported in $\Gamma$.
We  proved in  \cite{volkov2017reconstruction} the following theorem
\begin{thm} Let $\bg$ be in $\tilde{H}^{\f12}(\Gamma)^2$.
The problem
(\ref{uj1}-\ref{uj3andhalf})
%\bean
%\mu \Delta \bu+ (\lambda+\mu) \nabla \di \bu= 0  \mbox{ in } \RR^{3-} 
%\setminus \Gamma \label{uj1}, \\
% T_{\bev_3} \bu =0 \mbox{ on the surface } x_3=0 \label{uj2}, \\
% T_{\bn} \bu  \mbox{ is continuous across } \Gamma \label{uj3}, \\
% \aaa  \bu\bbb =\bg \mbox{ is a given jump across } \Gamma , \label{uj3andhalf}
%\bu (\bx) = O(\f{1}{|\bx|}),  \nabla \bu (\bx) = O(\f{1}{|\bx|^2}), \mbox{ uniformly as } 
%|\bx| \rightarrow \infty,
%\label{uj4}
%\eean
 has a unique solution in ${\cal V}$.
In addition, the solution $\bu$ satisfies the decay conditions
(\ref{uj4}).
%\bean
%\bu (\bx) = O(\f{1}{|\bx|^2}),  \nabla \bu (\bx) = O(\f{1}{|\bx|^3}), \mbox{ uniformly as } 
%|\bx| \rightarrow \infty.
%\label{uj4}
%\eean
\end{thm}

In this paper we will only consider forcing terms $\bg$ which are tangential to 
$\Gamma$. Physically, this reflects that the fault $\Gamma$ is not opening or starting to self intersect: only slip is allowed.
 We recall that if $\bg$ is continuous, the support of $\bg$, 
$\mbox{supp } \bg$, is equal to the closure
of the set of points in $\Gamma$ where $\bg$ is non zero; in general $\mbox{supp }  \bg$
is defined in the sense of distributions. 
%We will 
% assume that $\mbox{supp } \bg$ is exactly $\ov{\Gamma}$.
%This is not an overly restrictive assumption: if $\mbox{supp } \bg$ is strictly included in
% $\ov{\Gamma}$, if $\bg $ is smooth in $\Gamma$, and if the  interior ${\cal O}$, 
%relative to the topology
%induced in $\Gamma$ of 
%$\overline{ \{\bx \in \Gamma: \bg(\bx) \neq 0 \}}$ has a smooth boundary,
%we just re-define $\Gamma$ to be  ${\cal O}$.\\

\subsection{Fault inverse problem}
Can we determine both $\bg$ and $\Gamma$ from the data
 $\bu$ given only on the plane $x_3=0$?
Many investigators have studied
uniqueness and stability results for inverse boundary problems.
Earlier studies include papers such as Sylvester and Uhlmann's, 
\cite{sylvester1987global}, regarding the isotropic conductivity 
equation where it is proved that the knowledge of the Dirichlet
to Neumann boundary operator uniquely determines smooth conductivities.
In \cite{lee1989determining}, Lee and Uhlmann showed that
this is still true  in
  the  anisotropic case, up to a diffeomorphism.
On the subject of cracks, Friedman and Vogelius 
\cite{friedman1989determining} proved that, in dimension 2,
it suffices to apply two adequately chosen forcing terms on the boundary to uniquely determine  cracks in the framework of
the conductivity equation. 
The case of the two dimensional elasticity equation was considered 
by  Beretta et al. in \cite{beretta2008determination}. 
Stability results for linear cracks were derived; 
Beretta et al. proposed in \cite{beretta2010algorithm} a
MUSIC type algorithm for determining the position of these linear cracks
from boundary measurements. \\
 We note, however, that the case of interest in our present paper is substantially 
different for  
 two main reasons: first, the forcing term $\bg$ is given on the fault 
$\Gamma$, and second, our problem is three dimensional.
%involves 3D linear elasticity
%with constant Lam\'e coefficients. 
In \cite{volkov2017reconstruction}, we proved the following result:
\begin{thm}{\label{uniq1}}
Let $\Gamma_1$ and $\Gamma_2$ be two bounded open 
surfaces, with smooth boundary, such that each of them  is included in a
rectangle strictly contained in $\RR^{3-}$.
For $i$ in $\{ 1,  2\}$, assume that $\bu^i$ solves  (\ref{uj1}-\ref{uj4}) for 
$\Gamma_i$ in place of $\Gamma$ and $\bg^i$, a tangential field in 
$\tilde{H}^{\f12}(\Gamma_i)^2$, in place of $\bg$.
Assume that $\bg^i$ has full support in $\Gamma_i$, that is, 
$\mbox{supp } \bg_i = \ov{\Gamma_i}$.
Let $V$ be a non empty open subset in $\{x_3 =0\}$.
If $\bu^1$ and $\bu^2$ are equal in $V$, then
%\bean
% \Gamma_1 \cap \mbox{supp } \bg^1  =  \Gamma_2 \cap \mbox{supp } \bg^2
%\label{suppeq}
%\eean
$\Gamma_1 =  \Gamma_2$
and $\bg^1=\bg^2$.
\end{thm}
%\textbf{Remarks:} \\
%The notation $\mbox{supp }$ in (\ref{suppeq}) is for the support of $\bg^i$.
%A possible way of interpreting 
%the uniqueness result for the fault geometry (\ref{suppeq}) is to think of the two faults
%being equal on their "active" part.\\
There is a Green's tensor  $\bH$ such that 
		the  solution $\bu$ to problem (\ref{uj1}-\ref{uj3andhalf})  can also be written out as the convolution on $\Gamma$
\bean
  \int_\Gamma \bH(\bx, \by) \bg(\by) \, d \sigma (\by) \label{int formula},
\eean
The practical determination of this adequate half space Green's tensor $\bH$ was first studied
in  \cite{Okada} and later, more rigorously, in \cite{DV}.
Due to formula (\ref{int formula}) %and the  decay of the Green's tensor
%$H$ at infinity, 
we can define a continuous mapping ${\cal M}$ from tangential fields
$\bg $ in $\tilde{H}^{\f12}(\Gamma)^2$ to surface displacement fields $\bu(x_1, x_2, 0)$
in $L^2(V) $ where $\bu$ and $\bg$ are related
by (\ref{uj1}-\ref{uj4}). Theorem (\ref{uniq1}) asserts that 
this mapping is injective, so an inverse operator can be defined. % \cap C^\infty (\RR^2)$. 
It is well known, however, that such an operator ${\cal M}$ is compact, therefore
its inverse is unbounded. It is thus clear that any stable numerical method for reconstructing
$\bg$ from $\bu(x_1, x_2, 0)$ will have to use some regularization process.
In fact, in practice, our problem is even more challenging due to the fact that the geometry of the fault 
$\Gamma$ is also unknown. A numerical solution to determining $\Gamma$ and $\bg$ from 
$\bu(x_1, x_2, 0)$ will have to use a priori regularizing assumptions on $\bg$ 
and must be tested for robustness to noise.\\

\section{A functional for the  regularized reconstruction of  planar faults}
\label{Continuous regularized reconstruction} 
Let $R$ be a closed rectangle in the plane $x_3=0$.
Let $B$ be a set of $(a, b, d)$ such that the set
\bea 
 \{  (x_1, x_2, a x_1 + b x_2 +d): (x_1, x_2) \in R\}
\eea 
is included in the half-space $x_3 <0$.
We introduce the notations
\bea
m=(a,b,d), \\
\Gamma_m =\{  (x_1, x_2, a x_1 + b x_2 +d): (x_1, x_2) \in R\}.
\eea
We assume that $B$ is a closed and bounded subset of $\RR^3$. It follows that 
that 
\bean \label{pos dis}
\begin{array}{l}
\mbox{\sl the distance between $\Gamma_{m}$ and the plane $x_3=0$ is bounded below}
\\
\mbox{\sl by the same  positive constant for all $ m$ in $B$.}
\end{array}
\eean
In this section we assume that slips are supported in such sets $\Gamma_{m}$ (meaning that their supports are included in 
$\Gamma_m$, but they could be different from $\Gamma_m$).
We can then map all these fields into the rectangle $R$.
We thus obtain displacement vectors for $\bx$ in $V$ by the integral formula
\bean \label{new int}
\bu(\bx, \bg, m) = \int_R \bH_m(\bx, y_1, y_2)
 \bg (y_1, y_2) \sigma d y_1 d y_2, 
\eean
for any $\bg$ in $H^1_0 (R)$ and  $m$ in $B$, where $\sigma $ is the surface element on $\Gamma_m$
and  $\bH_m(\bx, y_1, y_2)$ is derived from 
the Green's tensor $\bH$ for $\by$ on $\Gamma_m$.
%Let $V$ be a non-empty, open subset of the plane
%$x_3 =0 $ and 
We now assume that $V$ is a bounded open subset of the plane $x_3 =0$ and for 
a fixed $\tilde{\bu}$ be in $L^2(V)$,
and  a fixed $m$ in $B$ we define the functional
\bean
 F_{m,C} (\bg)= \int_V  (\bu(\bx, \bg, m) - \tilde{\bu}(\bx ))'
{\cal C}^{-1}(\bx)  (\bu(\bx, \bg, m) - \tilde{\bu}(\bx ))
  d\bx
+ C \int_{R} |\nabla \bg|^2 ,
\label{F cont abd}
\eean
where ${\cal C} (\bx) $ is a diagonal positive definite 3 by 3 matrix 
for $\bx$ in $\ov{V}$, which is continuous in $\bx$, 
 and $C$ is a positive constant. In formula (\ref{F cont abd}) we intentionally
used ${\cal C}^{-1}$ rather than ${\cal C}$ because we will later view it as a covariance term.
 Define the operator
\bean 
A_{m} &:& H^1_0 (R) \ri L^2(V) \no \\
&&\bg \ri \int_R  \bH_m(\bx, y_1, y_2)
 \bg (y_1, y_2) \sigma d y_1 dy_2 .  \label{Aabd}
\eean
It is clear that $A_{m} $ is linear, continuous, and compact.
The functional $F_{m,C}$ can also be written as,
\bean
 F_{m,C} (\bg)= \|A_{m} \bg  - \tilde{\bu} \|^2_{L^2(V)} + C \| \bg \|^2_{H^1_0(R)},
\eean
where in $L^2(V)$ we use the norm
\bean \label{normV}
\| \bu \|_{L^2(V)} = 
(\int_V  \bu(\bx )'
{\cal C}^{-\f12}(\bx)  \bu(\bx) 
  d\bx)^{\f12},
\eean
and in $H^1_0(R)$ we use the norm
\bean \label{normR}
\| \bg \|_{H^1_0(R)} =
(\int_{R} |\nabla \bg|^2 )^{\f12}.
\eean
In the remainder of this paper, for the sake of simplifying notations,
both $\|   \q\|_{L^2(V)}$ and $\| \q  \|_{H^1_0(R)} $ will be abbreviated by $\| \q \|$;
context will eliminate any risk of confusion. 
\begin{prop}\label{Fmin}
For any fixed $m$ in $B$ and $C>0$, $F_{m,C} $ achieves a unique minimum $\bh_{m,C}$
in $H^1_0 (R)$.
\end{prop}
\textbf{Proof}:\\
The result holds thanks to classic Tikhonov regularization theory
(for example, see \cite{kress1989linear}, Theorem 16.4).

%On the subject of the widely used Tikhonov regularization technique, we recommend
%the textbook \cite{kress1989linear}, Chapter 16, for  a presentation 
%which is particularly  relevant to our study. 
For $\bh_{m,C}$  as in the statement of Proposition \ref{Fmin} we set,
	\bean \label{habdc def}
	 f_C(m) =  F_{m,C} (\bh_{m,C}).
	\eean
	\begin{prop} \label{fmin}
	$f_C$ is a Lipschitz continuous function on $B$. It achieves its minimum value on
	$B$.
\end{prop}
\textbf{Proof}:\\
We first note that the term $ \bH_m(\bx, y_1, y_2) $ and all its derivatives 
are uniformly bounded for $\bx$ in $V$, $(y_1, y_2)$ in $R$ and $m$ in $B$ thanks
to  (\ref{pos dis}). 
It follows that $ \bH_m(\bx, y_1, y_2) $ is Lipschitz continuous
in for $m$ in  $B$ with uniform Lipschitz constants for $(y_1, y_2)$ in $R$
and $\bx$ in $V$, so there is 
a positive constant $L$ such that
\bea
|\bH_m(\bx, y_1, y_2) - \bH_{m'}(\bx, y_1, y_2) | 
\leq L |m-m'|,
\eea
for any  $m$ and $m'$ in $B$,
for all $(y_1, y_2)$ in $R$, and all $\bx$ in $V$.
It follows that there is a constant ${\cal F}$  such that 
\bea
|F_{m,C} (\bh_{m,C}) - F_{m',C} (\bh_{m,C})| \leq {\cal F}|m-m'|,
\eea
for all $m$ and $m'$ in B.
By minimality for $\bh_{m',C}$, $ F_{m',C} (\bh_{m',C}) \leq F_{m',C} (\bh_{m,C})$,
so
\bea
F_{m',C} (\bh_{m',C}) \leq  F_{m,C} (\bh_{m,C}) + 
{\cal F}|m-m'|,
\eea
and given that we can switch the roles of $m$ and $m'$,
we found that
\bea
|F_{m',C} (\bh_{m',C}) -  F_{m,C} (\bh_{m,C})| \leq
{\cal F} |m-m'|.
\eea
Finally we just recall that $B$ is compact to claim that
$f_C$ achieves its minimum value.\\\\
	
The following theorem explains in what sense the argument of the minimum of the functional $F_{m,c}$ converges to the 	
	slip solving the fault inverse problem, and how the argument of the minimum of $f_C$ converges to the geometry parameter
	solving the fault inverse problem.
	
\begin{thm}\label{cv th} %948	
Assume that 
	$\tilde{\bu} = A_{\tilde{m}} \tilde{\bh}$ for
	some $\tilde{m}$ in $B$ and some 
	$\tilde{\bh}$ in  $H_0^1(R)$.
Let $C_n$ be a sequence of positive numbers converging to zero.
Let $m_n$ be any sequence in $B$ such that $f_{C_n}(m_n)$ minimizes 
 $f_{C_n}(m)$ for $m$ in $B$ and set
$f_{C_n}(m_n) = F_{m_n,C_n} (\bh_{m_n,C_n})$.
Then $m_n$  converges to $\tilde{m}$,
$\bh_{m_n,C_n}$ converges to $\tilde{\bh}$ in $H^1_0 (R)$,
and $A_{m_n} \bh_{m_n,C_n}$ converges to $\tilde{\bu}$ in $L^2(V)$.
\end{thm}
\textbf{Proof}:\\
We first note that 
\bean
\int_V | {\cal C}^{-\f12}(A_{m_n} \bh_{ m_n, C_n} - \tilde{\bu} )|^2  
+ C_n \int_{R} |\nabla \bh_{m_n,C_n}|^2 \no 
=f_{C_n}(m_n) \no \\%F_{m_n,C_n} (\bh_{m_n,C_n}) \no \\
\leq f_{C_n}(\tilde{m}) 
  \leq F_{\tilde{m},C_n} (\tilde{\bh})       
%\int_V | \bu(\bx, \tilde{\bh}, m_n) - \tilde{\bu}(\bx ) |^2  dx
 =C_n \int_{R} |\nabla \tilde{\bh}|^2. \label{contra ineq}
\eean 
Arguing by contradiction, assume that $m_n$ does not converge
to $\tilde{m}$.
After possibly extracting a subsequence, we may assume that 
$m_n$  converges to $m^*$
in $B$ with $m^* \neq \tilde{m}$.
By (\ref{contra ineq}) $\bh_{m_n,C_n}$ is bounded in $H^1_0(R)$: 
after possibly extracting a subsequence, we may assume
that $\bh_{m_n,C_n}$ is weakly convergent  to some $\bh^*$ in 
 $H^1_0(R)$.
Next we observe that $A_{m_n} $ is norm convergent to 
$A_{m^*}$ and we recall that $A_{m^*}$ is compact:
it follows that $A_{m_n} \bh_{m_n,C_n}$ converges strongly
to $A_{m^*} \bh^*$, thus we may take the limit as $n$ approaches
infinity of (\ref{contra ineq}) to find
\bea
\int_V | {\cal C}^{-\f12}(A_{ m^*} \bh^*- \tilde{\bu} )|^2   =0.
\eea
As $m^* \neq \tilde{m}$, this contradicts uniqueness
Theorem \ref{uniq1}. 
The same argument can be carried out to show that
$A_{m_n} \bh_{m_n,C_n}$ converges to $\tilde{\bu}$. \\
%We also argue by contradiction to show that $\bh_{m_n,C_n}$ converges to $\tilde{\bh}$:
%if that is not true  there is a positive $ \epsilon$ 
%such that for  a subsequence, still denoted by $\bh_{m_n,C_n}$ for convenience, 
%$\|\bh_{m_n,C_n} -  \tilde{\bh}\| \geq \epsilon$. 
We now show that $\bh_{m_n,C_n}$ converges to $\tilde{\bh}$ in $H^1_0 (R)$.
Since $A_{m_n} \bh_{m_n,C_n}$ converges to $\tilde{\bu}$, 
$A_{m_n}$ is  convergent to  $A_{\tilde{m}}$ in norm, 
and by
(\ref{contra ineq}) $\bh_{m_n,C_n}$ is bounded,
we can claim that 
$A_{\tilde{m}} \bh_{m_n,C_n}$ converges to $\tilde{\bu}$.
Let $\bv$ be in $L^2(\Gamma)$. The dot product in $L^2(V)$ and in $H^1_0(R)$
associated with the norms (\ref{normV}) and (\ref{normR}) will be denoted
by $(\q, \q)$.
As
\bea
(\bh_{m_n,C_n} - \tilde{\bh}, 
A_{\tilde{m}}^* \bv) = 
(A_{\tilde{m}} \bh_{m_n,C_n} - \tilde{\bu}, \bv)
\ri 0, 
\eea 
and $A_{\tilde{m}}$ is injective 
(due to theorem \ref{uniq1}, so the range of 
$A_{\tilde{m}}^*$ is dense), this shows 
that $\bh_{m_n,C_n}$ converges weakly to  $\tilde{\bh}$
in $H^1_0(\Gamma)$. 
To obtain strong convergence we recall that due to (\ref{contra ineq}),
$\| \bh_{m_n,C_n}\| \leq \| \tilde{\bh}\|$ and we write
\bean
\| \bh_{m_n,C_n} - \tilde{\bh}  \|^2  &=&
\| \bh_{m_n,C_n}\|^2 - 2 (\bh_{m_n,C_n},\tilde{\bh} ) +\|\tilde{\bh}  \|^2  \no \\
& \leq &2 (\tilde{\bh} -\bh_{m_n,C_n},\tilde{\bh} ), \label{good trick}
\eean
which tends to zero due to the weak convergence of $\bh_{m_n,C_n}$ to  $\tilde{\bh}$.

\section{A functional for the reconstruction of planar faults from a finite set of surface measurements} \label{finite recon}
%In this section we assume that $V$ is bounded.
For $j=1, .. ,N$, let $ P_j$ be points on the surface $x_3 =0 $ and
$\tilde{\bu}(P_j)$ be measured displacements at these points.
Let ${\cal F}_p$ be an increasing  sequence of finite-dimensional subspaces of
$H_0^1(R)$ such that $\bigcup_{p=1}^\infty {\cal F}_p$ is dense.
For $\bg$ in  ${\cal F}_p$ and $m$ in $B$, define the functional 
\bean
F^{disc}_{m,C} (\bg) = \sum_{j=1}^{N} C'(j,N) | {\cal C}^{-\f12}
((A_{m} \bg) - \tilde{\bu} ) (P_j)|^2 
+ C \int_R |\nabla \bg|^2 ,
\label{Ffunc}
\eean
where 
%$C'=((\sum_{j=1}^{N}  | \tilde{\bu} (P_j)  |^2))^{-1}$, $$
$A_{m} $ was defined in (\ref{Aabd}),  
$C>0$ is a constant. As to the constants $C'(j,N) $, simply put, they relate 
$F^{disc}_{m,C} $ to $F_{m,C} $ as $N$ tends to infinity. 
More precisely we assume that $\cal C$ is smooth and that 
for all positive integer $k$,  and 
for all $\varphi$ in $C^k(\overline{V})$, there is a constant 
$C(k)$ such that
\bean \label{quad rule}
|\int_V \varphi - \sum_{j=1}^N C'(j,N) \varphi(P_j) |
\leq C(k ) N^{-\beta} \sup_{V} \sum_{|l| \leq k}|D^l \varphi| ,
\eean
where $D^l \varphi $ is a  partial derivative of $\varphi$ with total order $l$
and $\beta$ is a positive integer depending on $k$.
We also assume that $C'(j,N) >0$ for all positive integer $N$ and all 
$j=1, .., N$.

	\begin{prop} \label{prop min}
The functional $F^{disc}_{m,C}$ achieves a unique minimum on ${\cal F}_p$.
\end{prop}
%	We don't know whether this minimum value is achieved at only one point.
%	Actually, for the value $N=1$ (only one measurement point), this is in general not true.
%	However, as $N$ grows large, we can state an interesting convergence result.
	\textbf{Proof}:\\
	This results again from Tikhonov regularization theory
(see \cite{kress1989linear}, Theorem 16.4).\\\\

According to Proposition \ref{prop min}, $F^{disc}_{m,C}$ achieves its minimum
	at some $\bh^{disc}_{m,C}$ in ${\cal F}_p$.
	We set 
	\bean
	 f^{disc}_C(m) =  F^{disc}_{m,C} (\bh^{disc}_{m,C}) \label{define bhdiscmc}.
	\eean
	\begin{prop}
	$f^{disc}_C$ is a Lipschitz continuous function on $B$ and achieves its minimum value on
	$B$.
\end{prop}

\textbf{Proof}:\\
The proof is similar to that of Proposition \ref{fmin}. \\

We now discuss the connection between the continuous and the discrete reconstruction functionals.
We will assume that the surface data $\tilde{\bu}$ is given by
$\tilde{\bu}= A_{\tilde{m}} \tilde{\bh}$, for some
 slip $\tilde{\bh}$ in $H^1_0(R)$ and $\tilde{m}$ in $B$.
Evidently, in the extreme case where the number of measurement points is $N=1$, we should expect 
no relation between $\bh^{disc}_{m,C}$ and $\tilde{\bh}$. In this section we want to analyze the convergence 
properties of $\bh^{disc}_{m,C}$ and the minimizer of $f^{disc}_C$ as the number of measurement points $N$ tends to infinity, 
$p$ tends to infinity, and $C$ tends to zero. Related proofs are rather intricate and technical, so we placed them in Appendix.

\begin{thm} \label{fdisc con}
Assume that 
	$\tilde{\bu} = A_{\tilde{m}} \tilde{\bh}$ for
	some $\tilde{m}$ in $B$ and some 
	$\tilde{\bh}$ in  $H_0^1(R)$.
Let $m^{disc}$ be such that
$$
f_C^{disc} (m^{disc}) = \min_B f^{disc}_C.
$$
Then for all $\eta >0$, there is an $N_0$ in $\NN$, a $p_0$ in $\NN$,
 and two positive constants $C_0$
 and $C_1$ such that if $N>N_0$, $p >p_0$, and $C_0 N^{-\beta}<C< C_1 $ then
\bean \label{disc geom conv}
| m^{disc} - \tilde{m}| \leq \eta.
\eean
\end{thm}
\textbf{Proof}:\\
The proof is given in Appendix.\\

\textbf{Remarks}:\\
Note that $m^{disc}$ depends implicitly on $N$, $C$, and $p$. \\
To interpret the condition $C_0 N^{-\beta}<C $ in Theorem \ref{fdisc con}, we recall that as $C$ tends to zero,
intuitively speaking, 
the functional $F^{disc}_{m,C}$ tends to an error (or misfit)
 calculated on the $N$ surface measurements 
$\tilde{\bu}(P_j)$.  Theorem \ref{fdisc con} states that a sufficient requirement for the reconstructed 
geometry parameter $m^{disc}$ to approach the real geometry parameter $\tilde{m}$
is for the regularization parameter $C$ to tend to zero and 
the subspace ${\cal F}_p$ to become large, all the while the number of measurement points $N$ tends to infinity 
with a rate such that   $C_0 N^{-\beta}<C $. Roughly speaking, this means that $N$ should not tend to infinity too slowly as $C$ tends to zero.\\

\section{Stochastic model} \label{Stochastic modeling}%247

\subsection{Model derivation}\label{Model derivation}
In our stochastic model we assume that the geometry parameter $m=(a,b,d)$ in $B$,  the 
 slip field $\bg$ in ${\cal F}_p$, and the measurements $\tilde{\bu} (P_j)$, are related by 
\bean \label{noise model}
 %m = \argmin_{ m \in B} \argmin_{\bg \in H^1_0(R)} F^{disc}_{m,C} (\bg)
(\tilde{\bu} (P_1), ..., \tilde{\bu} (P_N) )= 
(A_m \bg (P_1), ..., A_m \bg (P_N) )+ {\cal E},
\eean
where $A_m$ is given by (\ref{Aabd}), 
$m$ and $\bg$ are  now  random variables, and ${\cal E}$ in
$\RR^{3N}$ is additive noise, which 
is also assumed to be a random variable. 
We assume that ${\cal E}$ has a normal  probability density
$\rho_{noise}$ with mean zero and diagonal covariance matrix  such that %that for a measurement
 \bean  \label{rhonoise}
\rho_{noise}(\bv_1, ..., \bv_N) \propto \exp ( - \f12
 \sum_{j=1}^{N} C'(j,N) | {\cal C}^{-\f12}
\bv_j|^2 )
\eean 
Accordingly, the  probability density of the measurement
$\tilde{\bu}_{meas}$ knowing the geometry parameter $m$ and the slip field $\bg$ 
is 
\bean \label{umeas}
\rho(\tilde{\bu}_{meas}  | m, \bg ) \propto  
\exp ( -\f12 \sum_{j=1}^{N} C'(j,N) 
| {\cal C}^{-\f12}
(A_m \bg- \tilde{\bu}_{meas})(P_j) |^2 ).
\eean
Next, we assume that the random variables $m$ in $B$ and $\bg$ in 
${\cal F}_p$  are independent.
The prior distribution of $m$, $\rho_{prior}$ is assumed to be uninformative,
that is, $\rho_{prior}(m) \propto 1_B (m)$ and 
the prior distribution of  $\bg$ is 
assumed to be Gaussian with mean zero and given by
\bean \label{gprior}
\rho_{{\cal F}_p}(\bg) \propto \exp( -\f12 C \int_R |\nabla\bg|^2).
\eean
Applying the Bayesian theorem and independence of the priors of $m$ and $\bg$, 
we write
\bean
\rho( m, \bg| \tilde{\bu}_{meas}  ) \propto 
\ds \rho(\tilde{\bu}_{meas}  | m , \bg) \rho_{{\cal F}_p}(\bg) \rho_{prior} (m)
.
\label{Bayes}
\eean
In this work we are only interested in recovering the posterior probability density of $m$,
so we integrate formula (\ref{Bayes}) in $\bg$ over ${\cal F}_p$ to obtain the marginal posterior 
density of $m$. It turns out that this can be done explicitly
thanks to the Gaussian formulation.
Introducing the $3N$ by $3N$ diagonal  matrix
 ${\cal D}$  such that
\bean \label{calD def}
{\cal C}^{-\f12} (C'(1,N)^{\f12}\bu (P_1), ..., C'(N,N)^{\f12}\bu(P_N) ) = {\cal D}
(\bu (P_1), ..., \bu(P_N) ),
\eean
we can state,
\begin{prop} \label{int g prop}
Integrating in $\bg$ over ${\cal F}_p$ the right hand side of formula
(\ref{Bayes}), we find 
\bean  \label{m dist}
\rho( m| \tilde{\bu}_{meas}  ) \propto
\exp( -\f12 F^{disc}_{m,C} (\bh^{disc}_{m,C})) 
\f{\rho_{prior} (m) }{\s{\det ((2\pi)^{-1} (A_m'{\cal D}^2A_m + C I_q))}},
\eean
where $F^{disc}_{m,C}$ is given by (\ref{Ffunc}),
$\bh_{m,C}^{disc}$ is defined by (\ref{define bhdiscmc}), $I_q$ is the identity operator of the 
$q$ dimensional
subspace ${\cal F}_p$, $A_m$, initially defined by (\ref{Aabd}), is here restricted to a linear operator
from ${\cal F}_p$ to $\RR^{3N}$.
\end{prop}
\textbf{Proof:}\\
The proof is given in Appendix.

\subsection{Proof of convergence of stochastic model to the unique solution of the deterministic inverse problem
as covariance tends to zero}\label{convergence of stochastic model}
Investigators have examined the limit probability laws 
for inverse problems as the number of measurement points grows large
\cite{ammari2011imaging} %ammari garnier 
or as the random fluctuations of the media are small
\cite{borcea2003theory} % Borcea Tsogka equation 2.6
but these studies pertained to source location in media with propagating waves, so the 
connection to our framework is non trivial and a detailed mathematical analysis is likely to
be involved. We can still state and prove the following result regarding the pointwise
convergence of the posterior distribution $\rho( m| \tilde{\bu}_{meas}  )$ 
for small noise covariance. The proof relies on estimates shown in the previous two sections.
Formally,  the limit case where the covariance  defined by ${\cal C}$ tends to zero
can be interpreted in light of Theorem \ref{fdisc con}. 
Suppose that  the measurements $\tilde{\bu}_{meas} (P_j)$ were produced by a slip on a fault whose
geometry was given by $\tilde{m}$. %plus some noise, a realization of ${\cal E}$.
As ${\cal C}$  tends to zero,  $\rho( m| \tilde{\bu}_{meas}  )$  tends formally to the Dirac measure centered
at $\tilde{m}$, so the probability density  $\rho( m| \tilde{\bu}_{meas} )$
will achieve its maximum arbitrarily close to $\tilde{m}$.
%, for all large enough $N$ and 
%small enough $C$. 
More precisely, we now set
 \bean \label{posterior}
%\f{\ds \rho(\tilde{\bu}_{meas}  | m ) \rho_{prior} (m)}{ \rho(\tilde{\bu}_{meas} )} 
\rho_\tau (m | \tilde{\bu}_{meas} )=  
{\cal I}_\tau \exp( -\f{\tau}{2} F^{disc}_{m,C} (\bh^{disc}_{m,C})) \f{ 1_B(m)}{\s{\det ((2\pi)^{-1} 
\tau (A_m'{\cal D}^2A_m + C I_q))}},
\eean
where $\tau>0$ is a constant that will tend to infinity, $1_B$ is the indicator function of $B$,
and ${\cal I}_\tau $ is a normalizing constant. Note  the new surface covariance is 
 ${\tau^{-1}\cal C}$, so letting $\tau \ri \infty$ will
make this rescaled covariance tend to zero. Observe also that in this formulation, 
since both ${\cal C}^{-1}$ and $C$ are both rescaled by $\tau$, $\bh^{disc}_{m,C}$ is independent
of $\tau$.

\begin{prop}
Suppose that the probability density of $m$ knowing surface measurements is given by 
(\ref{posterior}) where $\tilde{\bu}$ in the definition (\ref{Ffunc}) of  $F^{disc}_{m,C}$ is given by
$\tilde{\bu}= A_{\tilde{m}} \tilde{\bh}$ for some $\tilde{m}$ in $B$ and $\tilde{\bh}$ in $H^1_0(R)$
and $\bh_{m,C}^{disc}$ is as in (\ref{define bhdiscmc}).
%measurements $\tilde{\bu}_{meas} $ are given by 
%$A_{\tilde{m}} \tilde{\bh}(P_j)$ for some $\tilde{m}$ in $B$ and $\tilde{\bh}$ in $H^1_0(R)$ plus a realization
%of the noise $\cal E$.
Let $m_0$ be in $B$ such that  $m_0 \neq \tilde{m}$.
 Then for all large enough $N$ and $ p$, and small enough $C$, the posterior distribution $\rho_\tau( m| \tilde{\bu}_{meas}  )$  
evaluated at $m_0$ converges to zero 
 as $\tau \ri \infty$. Additionally, 
this convergence is uniform in  $m_0$,  as long as $m_0$ remains bounded away
from $\tilde{m}$.
\end{prop}

\textbf{Proof}:\\
Fix $m_0$ different from $\tilde{m}$. 
It is clear that the exponential and the fraction
in formula (\ref{posterior}) converges to zero as $\tau$ tends to infinity, so the crux of the proof is to account for the normalizing constant ${\cal I}_\tau$.
According to Theorem \ref{fdisc con},
% with $s^{disc}_C$ in place of $f^{disc}_C$,
%there is an $N_0$ in $\NN$ and a positive $C_0$
%such that if $N>N_0$ and $0<C<C_0$ then 
for a fixed $N$, $p$, and $C$, provided $N$ and $p$ are large enough and $C$ is small enough, 
the minimum of $f_C^{disc}$ will occur
in the ball with center $\tilde{m}$ and radius $\f{|\tilde{m} - m_0 |}{2}$.
We now fix such an $N$, $p$,  and $C$.
Let $ m '$ be a point where  $f_C^{disc}$  will achieve its minimum in this ball.
Let $\bh_{m_0}$ be the element in ${\cal F}_p$ where $F_{m_0,C}^{disc}$ achieves its minimum
and $\bh_m'$ be the element in ${\cal F}_p$ where $F_{m',C}^{disc}$ achieves its minimum.
Set $\gamma = F_{m_0,C}^{disc}(\bh_{m_0}) - F_{m',C}^{disc} (\bh_m')$. We must have $\gamma >0$ since
$m_0$ is not in the ball with center $\tilde{m}$ and radius $\f{|\tilde{m} - m_0 |}{2}$.
But by continuity, there is a positive $\alpha$ such that if $|m - m'| \leq \alpha$,
\bean
F_{m_0,C}^{disc}(\bh_{m_0}) - F_{m,C}^{disc} (\bh_m') \geq \f{\gamma}{2}.
\eean
Let $\bh_m$ be the element in $H^1_0(R)$ where $F_{m,C}^{disc}$ achieves its minimum.
Necessarily, if $|m - m'| \leq \alpha$, 
\bean
F_{m_0,C}^{disc}(\bh_{m_0}) - F_{m,C}^{disc} (\bh_m) \geq \f{\gamma}{2}.
\eean
We can now estimate the normalizing constant ${\cal I}_\tau$.
Now letting $\alpha$ be the maximum of the largest eigenvalue of $A_m'{\cal D}^2A_m $ for 
$m$ in $B$, 
\bea
\f{ 1}{\s{\det ((2\pi)^{-1} 
\tau (A_m'{\cal D}^2A_m + C I_q))}} \geq \f{1}{((2 \pi)^{-1} \tau (C+\alpha)^{\f{q}{2}}},
\eea
where $q=\dim {\cal F}_p$. Since we have 
\bea
\int_B \exp (  - \f{\tau}{2}F_{m,C}^{disc}(\bh_{m})  ) \f{ 1}{\s{\det ((2\pi)^{-1} 
\tau (A_m'{\cal D}^2A_m + C I_q))}}  dm \\
\geq \f{1}{((2 \pi)^{-1} \tau (C+\alpha))^{\f{q}{2}}}
  \int_{B \cap \{ |m' - m| \leq \alpha \}} 
\exp ( \tau \f{\gamma}{4}  - \f{\tau}{2}F_{m_0,C}^{disc}(\bh_{m_0})  ) dm,
\eea
and
we have,
\bea
{\cal I}_\tau = O( \exp( -\tau \f{\gamma}{4}  + \f{\tau}{2}F_{m_0,C}^{disc}(\bh_{m_0}))
((2 \pi)^{-1} \tau (C+\alpha))^{\f{q}{2}}).
\eea
It follows that
\bea
\rho_\tau( m_0| \tilde{\bu}_{meas}  ) =  O( \exp( -\tau \f{\gamma}{4} ) 
  (1+\f{\alpha}{C})^{\f{q}{2}}),
\eea
which tends to zero as $\tau \ri \infty$.
We note that this convergence is uniform in $m_0$,
 as long as $m_0$ remains bounded away
from $\tilde{m}$: this is due to Theorem \ref{fdisc con}.

\subsection{Discrete problem and size of computation}
As mentioned earlier, $H^1_0(R)$ is approximated by the  $q$ -dimensional vector space over $\RR$, ${\cal F}_p$. 
The slip field $\bg$ will be approximated by $g^{(p)}$ and the term $\nabla \bg $ in (\ref{Ffunc})
will be approximated by $D g^{(p)}$ for $\p_{y_1} \bg$
and $E g^{(p)}$ for $\p_{y_2} \bg$ where $D$ and $E$ are two invertible $q$
by $q$ matrices. 
The term $(A_{m} \bg) (P_j)$, $j=1, .., N$, in (\ref{Ffunc}) is approximated
by $A g^{(p)}$, where $A$ is a $3N \times q$ matrix (the reader should bear in mind that
this matrix depends on $m$).
The singular values of the continuous operator $A_{m}$ are known
to decrease very fast, \cite{GohbergKrein}.
Accordingly, the singular values of $A$ decrease fast too. 
The discrete equivalent of minimizing (\ref{Ffunc}) is to minimize
\bean
 \| {\cal D} (A g^{(p)} - u^{(3N)})\|^2 + C (\| D g^{(p)}  \|^2 + \| E g^{(p)}  \|^2),\label{disc pb}
\eean
where we use the usual Euclidean norms on $\RR^{3N}$ and on $\RR^q$
and $u^{(3N)}$ in $\RR^{3N}$ stands for the data 
$  \tilde{\bu}_{meas}(P_j)$, $j=1, .., N$.
If $C>0$ is known, this amounts to solving the following linear equation
\bean
A'{\cal D}^2 A g^{(p)} + C (D' D  + E' E )g^{(p)} = A'{\cal D}^2 u^{(3N)} \label{disc opt}.
\eean
In practice this inverse problem is vastly underdetermined since
$3N << q$. Even if  $3N$ was greater or equal than $q$, it would not be 
possible to set $C=0$ since the singular values of $A$ decay very fast. 
We  set a grid of points in $B$
\bea
(a_{i_1},b_{i_2},d_{i_3}), (i_1,i_2,i_3) \in I.
\eea  
Thus all together, if $C$ is known the $q \times q$ linear system (\ref{disc opt})
has to be solved $|I|$ times. 
Here we recall $R$ is two dimensional;  in our practical calculations $q$ was $50^2$ or $100^2$, while
$|I|$ was about $50^3$. Finally, an appropriate value for $C$ had to be computed by iterations, so
for each value of $(i_1,i_2,i_3)$ the  linear system (\ref{disc opt}) had to be solved multiple times (the number of iterations  depended on $(i_1,i_2,i_3)$, it was on a range from 
5 to 50, so all together linear system  (\ref{disc opt}) had to be solved
about  $10^6$ to $10^7$ times).
% This accounts for the very large size of the computation despite the fact that for each $m$,
%$\tilde{\bh}$ was assumed to be deterministic.
In fact our calculation was only possible to perform thanks to the use of adequate  linear 
algebra techniques and a parallel multi core implementation which are described in details in
subsequent paragraphs.

\subsection{Algorithm for selecting the regularizing constant $C$}
%This kind of question has been a long-standing conundrum 
%in the field of inverse problems. 
Some general suggestions  for selecting the regularizing constant $C$
can be found in the literature,
\cite{Golub1979, Hansen1992, Kilmer2001}.
We note, however, that some well known methods (for example, truncated SVD) are inapplicable
 since we are not within the framework of $L^2$ regularization. 
%,
%and we search for constant sign solutions. 
The additional challenge in our case is that we have as many matrices $A$ as 
different points $(a_{i_1},b_{i_2},d_{i_3})$ for the chosen grid of $B$.
Our method  for selecting a constant $C$ uses the following lemma.
\begin{lem} \label{C exists}
Let ${\cal D}v^{(3N)}$ be the orthogonal projection of ${\cal D} u^{(3N)} $
on the range of ${\cal D}A$. Assume that $v^{(3N)}$ is non-zero.
Let $g^{(p)}$ solve (\ref{disc opt}).
Then $\| {\cal D}(Ag^{(p)} - u^{(3N)})\|$ is a continuous function of $C$ in $(0, \infty)$ with range
$(\|{\cal D}(u^{(3N)}- v^{(3N)})\|, \|{\cal D} u^{(3N)}\|)$.
\end{lem}
\textbf{Proof}:\\
The proof is given in Appendix.

\begin{definition} \label{defC} %1306
With the same notations as in Lemma  \ref{C exists}, set $\mathrm{{\bf Err}}$
be a number in $(0, \| {\cal D} u^{3N} \|)$. For every $(i_1,i_2,i_3) \in I$ we set
$C(i_1,i_2,i_3) =0$, if $\mathrm{{\bf Err}} \leq \|{\cal D} (u^{(3N)}- v^{(3N)})\| $,
otherwise we set
\bea
C(i_1,i_2,i_3) = \sup \{C>0: \| {\cal D} (Ag^{(p)} - u^{(3N)})\| \leq \mathrm{{\bf Err}}  \}.
\eea
\end{definition}

Solving for $C(i_1,i_2,i_3)$ for $(i_1,i_2,i_3) \in I$
is expensive since in practice $I$ has a large cardinality,
and for each $(i_1,i_2,i_3)$ in $I$
a non-linear inequality has to be solved where at each iterative step a large linear system
has to be solved.

We employed  several techniques for cutting down computational time.
\begin{enumerate}
\item The Green's tensor for half space linear elasticity $\bH(\bx,\by)$ is notoriously expensive to compute. It is possible, however, to use a simplified form since we know that $x_3 =0$. 
Formulas that reduce the number of operations  were given in \cite{DV}.
Additional savings in computational time were achieved for the fault inverse problem as 
we computed values $\bH(\bx,\by)$ by passing  a single, large array.%, for each  $(i_1,i_2,i_3)$
\item Woodbury's formula is remarkably helpful for a fast computation of \\
$(A'{\cal D}^2A + C (D' D  + E' E ))^{-1}A'{\cal D}^2u^{(3N)}$: this is because $(D' D  + E' E)^{-1}$ can be
 pre-computed
and stored, and $3N << q$, so computing  the SVD of ${\cal D} A$ is cheap, and $A'{\cal D}^2A$ has low rank.
\item Finally, since the set of indices $I$ is typically large, it is greatly beneficial
to use a multi-core parallel implementation. After pre-computing a few variables,
 matrices, and tables, computing all the constants $C(i_1,i_2,i_3)$ can be done in parallel.
\end{enumerate}
This led us to the following algorithm.
\begin{center}
\begin{tabular}{ |l |}
\hline
\textbf{Algorithm for computing $C(i_1,i_2,i_3)$, $(i_1,i_2,i_3) \in I$ } \\
\hline
Compute and save $(D' D  + E' E)^{-1}$\\
For each $(i_1,i_2,i_3) \in I$\\
\phantom{XX} Compute $A$\\
\phantom{XX} Compute the SVD of ${\cal D}A$\\
\phantom{XX} if $\|{\cal D}(u^{(3N)}- v^{(3N)})\| < \mathrm{{\bf Err}} $\\
\phantom{XX} \phantom{XX} use a non-linear iterative solver to find $C(i_1,i_2,i_3)$\\
\phantom{XX} \phantom{XX} \% \textsl{at each iteration use Woodbury's formula
for computing $(A'{\cal D}^2A + C (D' D  + E' E ))^{-1}A'{\cal D}^2u^{(3N)}$}\\
\phantom{XX} else \\
\phantom{XX} \phantom{XX} $C(i_1,i_2,i_3) =0$\\
\phantom{XX} end  \\
end\\
\hline
\end{tabular}
\end{center}
We are now able to define a uniform regularization constant $C$ by setting
\bean \label{choice C}
  {\bf C} = \max_{(i_1,i_2,i_3) \in I}{ C(i_1,i_2,i_3)}.
\eean
To justify this choice, first we
 note that for each $(i_1,i_2,i_3)$ in $I$, intuitively speaking, $C(i_1,i_2,i_3)$ will give 
rise to the most regular (or the least energy) pseudo-solution to the equation
$A g^{(p)} = u^{(3N)}$ with error tolerance $\mathrm{{\bf Err}}$.
This is in line with previous studies of de-stabilization of faults 
\cite{dascalu2000fault, ionescu2008earth, volkov2010eigenvalue}
and general minimum energy principles in physics.
Next, picking the same $ {\bf C}$ for all $(i_1,i_2,i_3) \in I$
will make it possible to find the geometry 
$m$ that will best replicate 
the surface measurements $\tilde{\bu}_{meas} (P_j)$ with a uniform 
control of
energy exerted by the slip.

%\subsection{Algorithm for computing the probability density
%$\rho(\bigwedge_{j=1}^N \bv (P_j) | m )$}
\subsection{Algorithm for computing the probability density
$\rho(m|u^{(3N)})$ }
The probability density of $\rho(m|u^{(3N)})$ is given by
\bean \label{gauss prob dens}
 {\cal I} \exp(-\f12 \| {\cal D} (A g^{(p)} - u^{(3N)})\|^2
- \f12 {\bf C} (\| D g^{(p)} \|^2 + \| E g^{(p)} \|^2)) \no \\
\f{\rho_{prior} (m) }{\s{\det ((2\pi)^{-1} (A'{\cal D}^2A + {\bf C}  I_q))}},
\eean
where ${\cal I}$ is a normalizing constant (which is unknown at the beginning of the calculation) and 
$$
g^{(p)}= (A'{\cal D}^2A + {\bf C} (D' D  + E' E ))^{-1} A'{\cal D}^2 u^{(3N)}. 
$$
Again, since the cardinality of $I$ is rather large, we are careful to pre-compute and store 
adequate terms. This led us to the following algorithm.
\begin{center}
\begin{tabular}{ |l |}
\hline
\textbf{Algorithm for computing the probability density
$\ds \rho(u^{(3N)} | m(i_1,i_2,i_3) ) $,
$(i_1,i_2,i_3) \in I$} \\
\hline
Compute and save $(D' D  + E' E)^{-1}$ \\%and $ {\cal C}^{-1}$\\
For each $(i_1,i_2,i_3) \in I$\\
\phantom{XX} Compute $A$\\
\phantom{XX} Compute the SVD of ${\cal D}A$\\
\phantom{XX} Use  the singular values  of ${\cal D}A$ to compute $\det ((2\pi)^{-1} (A'{\cal D}^2A + 
{\bf C}  I_q))$\\
\phantom{XX} Solve $A'{\cal D}^2A g^{(p)} + {\bf C} (D' D  + E' E )g^{(p)} = A'{\cal D}^2u^{(3N)} $ \\
%\phantom{XX} Evaluate  $v^{(3N)}= A g^{(p)}$\\
\phantom{XX} Evaluate $\ds \rho( m|  u^{(3N)})/{\cal I}$ by using
 formula (\ref{gauss prob dens}) \\
end\\
\hline
\end{tabular}
\end{center}
After applying the algorithm, there only remains to compute the normalizing constant
${\cal I}$. 
Since 
$A_{m} $ is a smooth 
function of $m$  it follows that $g^{(p)}$ is in turn smooth in $m$,
and so is $\ds \rho( m|  u^{(3N)})/{\cal I}$. 
The exponential in formula  (\ref{gauss prob dens}) shows that 
$\ds \rho( m|  u^{(3N)})/{\cal I}$ is also rapidly decaying way from its maximum,
so $\int \ds \rho( m|  u^{(3N)})/{\cal I} dm$ can be efficiently evaluated 
by the three dimensional trapezoidal rule.

%Even though computing 
% We  note that the integral
%formula (\ref{Bayes}) is evaluated numerically in our calculation.
%The accuracy depends on how fine the grid of $B$ $(a_{i_1},b_{i_2},d_{i_3})$ for
%$(i_1,i_2,i_3) \in I$.
%In practice if $B$ is adequately selected, the probability density 
%$\rho(\bigwedge_{j=1}^N \tilde{\bu}(P_j) | m) $ is smooth and decays very fast as we move
%toward the boundary of $B$, therefore using a straightforward
%regular grid of $B$ will lead to an accurate evaluation of the integral term
%in (\ref{Bayes}).

\section{Numerical results}\label{Numerical results} %268
We present in this section numerical results for the recovery of the fault geometry
parameter $m$
% and slip field
from surface measurements.
Ultimately, we will show results pertaining to the specific case of the 2007 
Slow Slip Event (SSE) which
occurred in  the Guerrero region of Central Mexico. 
\begin{figure}[ht]
\begin{center}
\includegraphics[scale=.5]{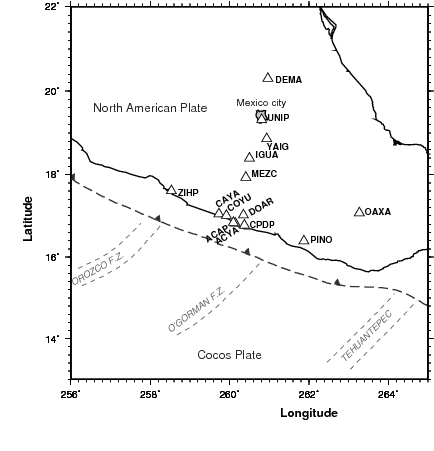}
\caption{The Guerrero gap region of Mexico. The subduction zone studied in this paper
meets the sea floor of the Pacific ocean along a nearly linear course called
the Middle American Trench: it appears on this figure 
as a dashed line. The large triangles mark
the locations of  GPS stations that were used to record the 2006 Guerrero SSE.
%They are named ACAP, ACYA, CAYA, COYU, CPDP, 
%DEMA, DOAR, IGUA, MEZC, OAXA, 
%PINO, UNIP, YAIG, ZIHP
} \label{geography}
\end{center}
\end{figure}
Figure \ref{geography} shows locations where surface displacements have been
 recorded. Over time,
some recording stations have closed while others have opened: it is therefore preferable not to use
all these stations in the training step. Specifically, 
we will use  measurements from the stations
ACAP, ACYA, CAYA, COYU, CPDP, DEMA, DOAR, IGUA, MEZC, UNIP, YAIG.
%see Figure  \ref{geography}.
In effect, these will be the points $P_j$ in our computations, for $j=1, .., N$, and
$N$ is equal to 11.
We use a rectangular system of coordinates $x_1, x_2, x_3$ centered
at ACAP: the $x_1$ direction runs West-East, the $x_2$ direction runs
South-North, and the $x_3$ direction runs down-up.
 In effect, this assumes that the Earth is locally flat. 
 Units for distances will be kilometers. Local geography is ignored,
 so $x_3=0$ at each of these 11 stations.
The medium  Lam\'e coefficients $\lambda$ and $\mu$ 
will be set to 1, which results in a Poisson ratio 
0.25, a commonly agreed upon value for Earth's rocks.
We refer to \cite{volkov2017determining, volkov2017reconstruction}
for an account of how raw displacement data was collected day after day.
The data was then completed and smoothed, as explained in \cite{volkov2017determining}.
The error bar on the data can be estimated by comparing the smoothed data to the raw data.
Here we have  to emphasize that finding the most optimal  and  accurate estimates
of the average and the standard deviation of displacement fields  is beyond the scope of our work.
However, satisfactory estimates are easy to find and provide a good starting point 
for addressing the stochastic fault inverse problem.
The effective maximum of $| \tilde{\bu}(P_j)|$ is about 100 mm.
The standard deviation on measurements of horizontal displacements  can be estimated to
 .8 mm, and 2 mm for the 
 standard deviation 
 of vertical displacements.
We will show in this section three test cases before covering
the real world case. In the test cases, the surface points $P_j$ will be the same as the ones used in the real world data case. We simulated data and added gaussian noise
with same covariance as the one estimated in the real world data case.
In the test cases we made sure to set faults  with depths that are consistent with what 
geophyscists expect to find in that region (in general, these depths are not deeper than 80 km, due
to the thickness of Earth's crust) and to produce surface displacements with the same order
of magnitude as those observed for the 2007 Guerrero SSE.\\
In each case we set the center of the rectangle $R$ to be the average
of the coordinates of $P_j$ weighted by $|\tilde{\bu} (P_j)|$.
The lengths of the sides of the rectangle $R$ can be set by first examining a large area which
includes all the $P_j$s  and then re-focusing it from a reconstructed $\bh$.
Alternatively, the size of rectangle $R$ can be estimated by applying the quasi constant slip method
presented in section 3.1 of \cite{volkov2017determining}.

\subsection{First test case}
\begin{figure}[ht] 
\begin{center}
\includegraphics[scale=.4]{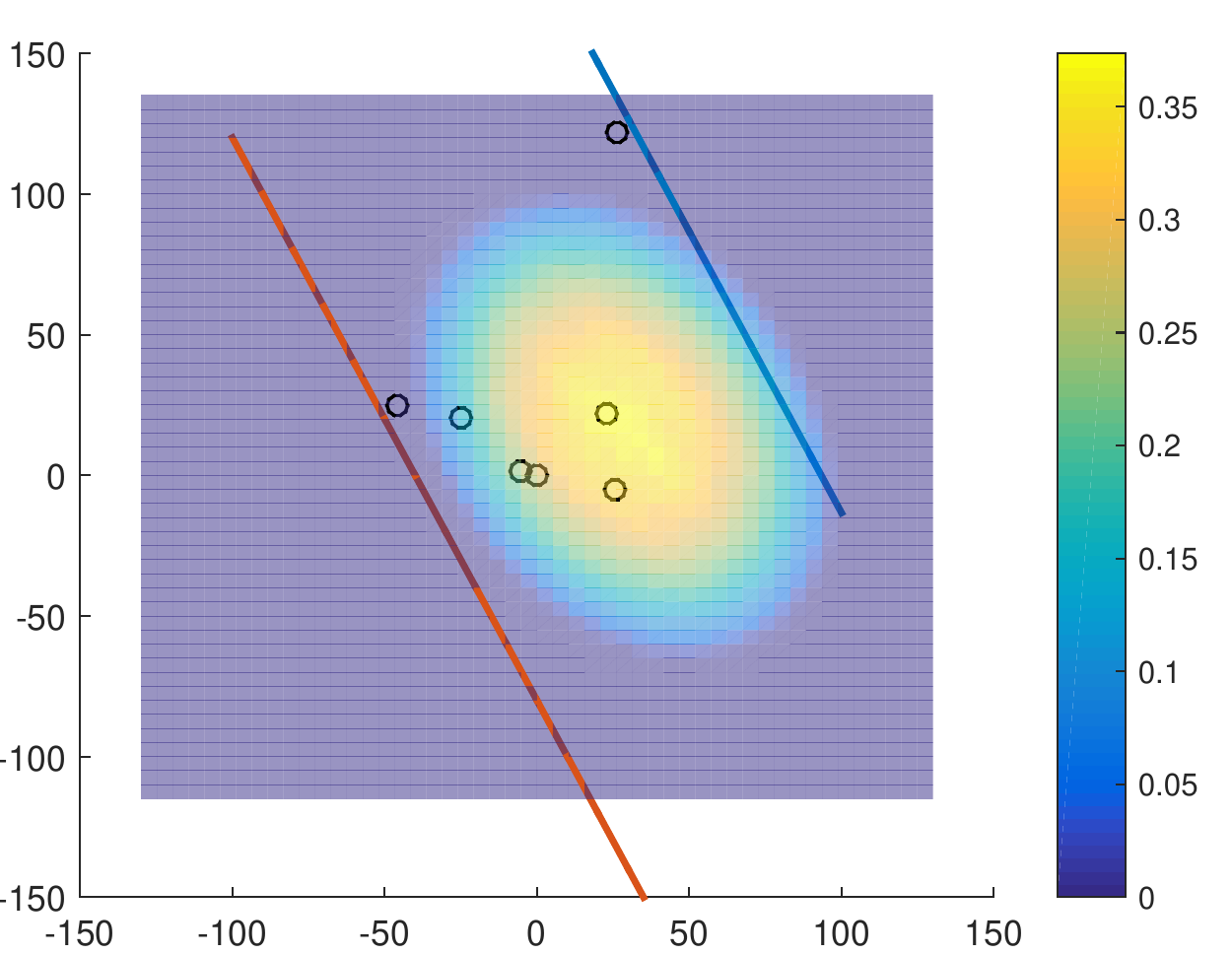}
\includegraphics[scale=.4]{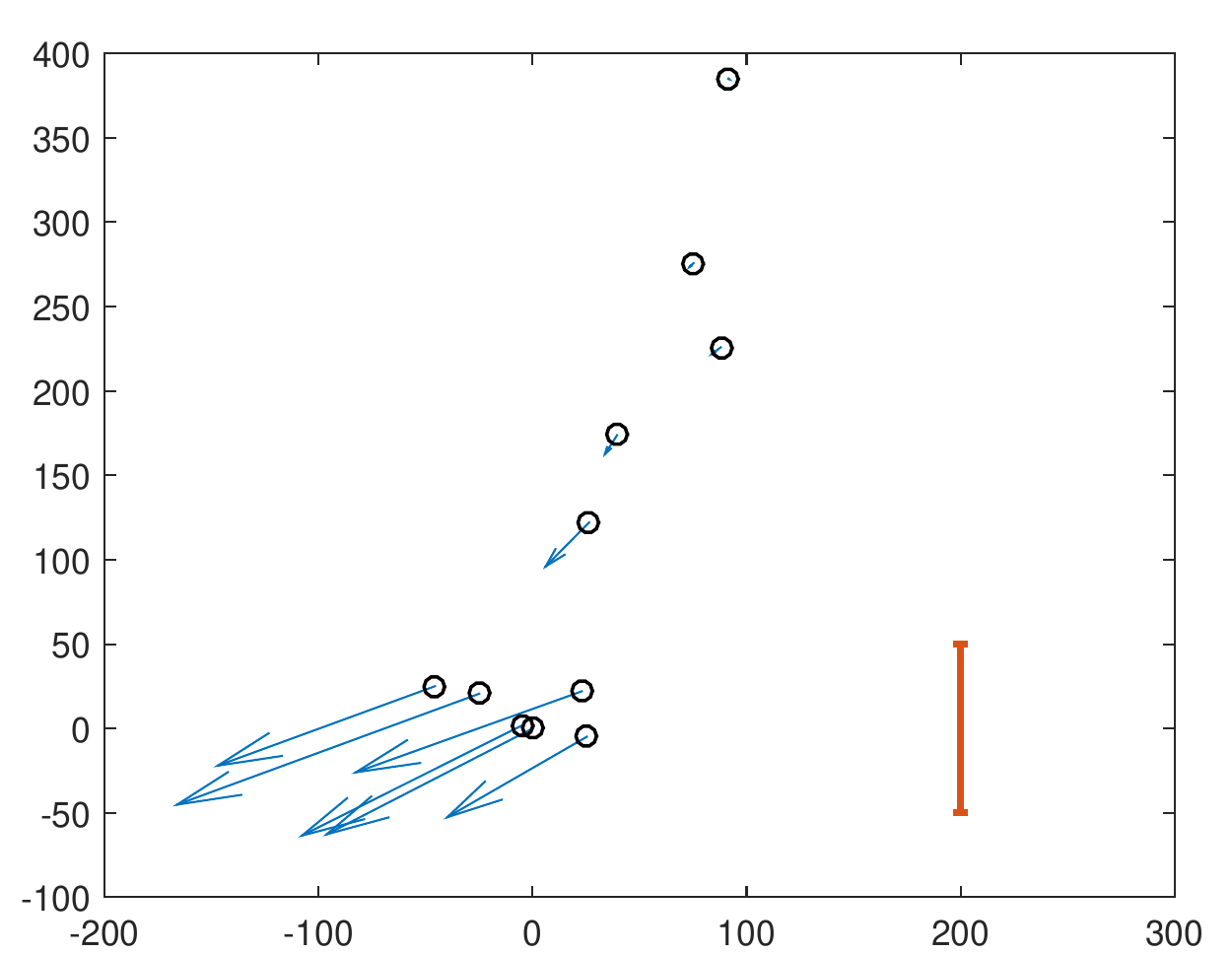}
\includegraphics[scale=.4]{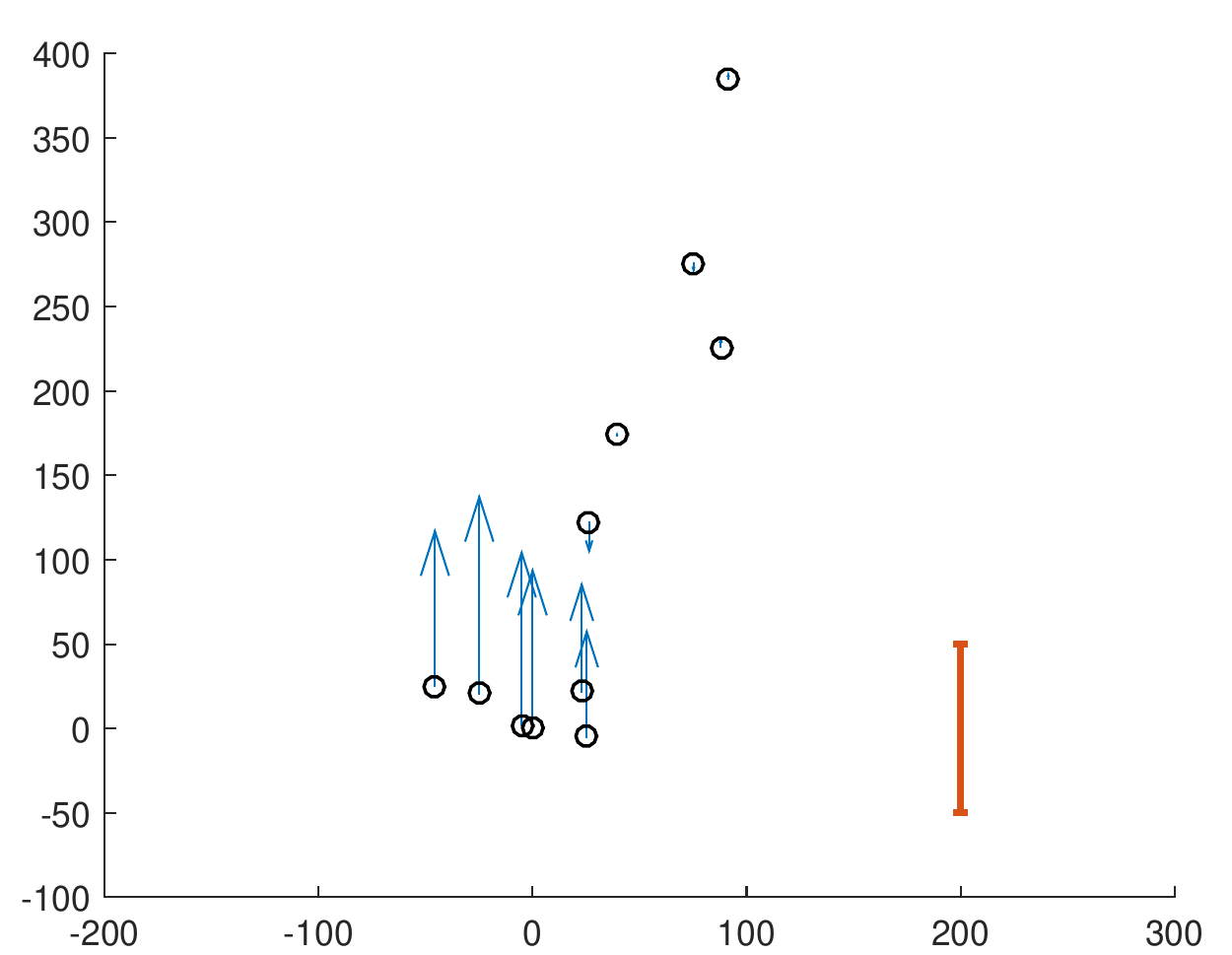}
\caption{Test case 1. Top left: the fault $\Gamma$ and
the slip field $\tilde{\bh}$. The red line is $x_3=-2$, the blue line
 is $x_3=-40$, on the plane $x_3=ax_1 + bx_2 + d$.
The circles stand for the surface measurement points $P_j$. They appear
on the map in Figure \ref{geography}. Units for surface distances are  kilometers. 
Color bar shows $|\tilde{\bh}|$, in meters. $\tilde{\bh}$ points in the direction of steepest 
ascent. The next two panels show the resulting surface displacements
at the $P_j$s. The red line segment indicates  the scale: 100 mm.} \label{geo15}
%Top right: reconstruction. Bottom two graphs show horizontal and %vertical displacements,
%data $\bu(P_j)$ in red and reconstructed surface field
%$\tilde{\bu}(P_j)$ in blue  } \label{test2fig}
\end{center}
\end{figure}

\begin{figure}[ht] 
\begin{center}
\includegraphics[scale=.3]{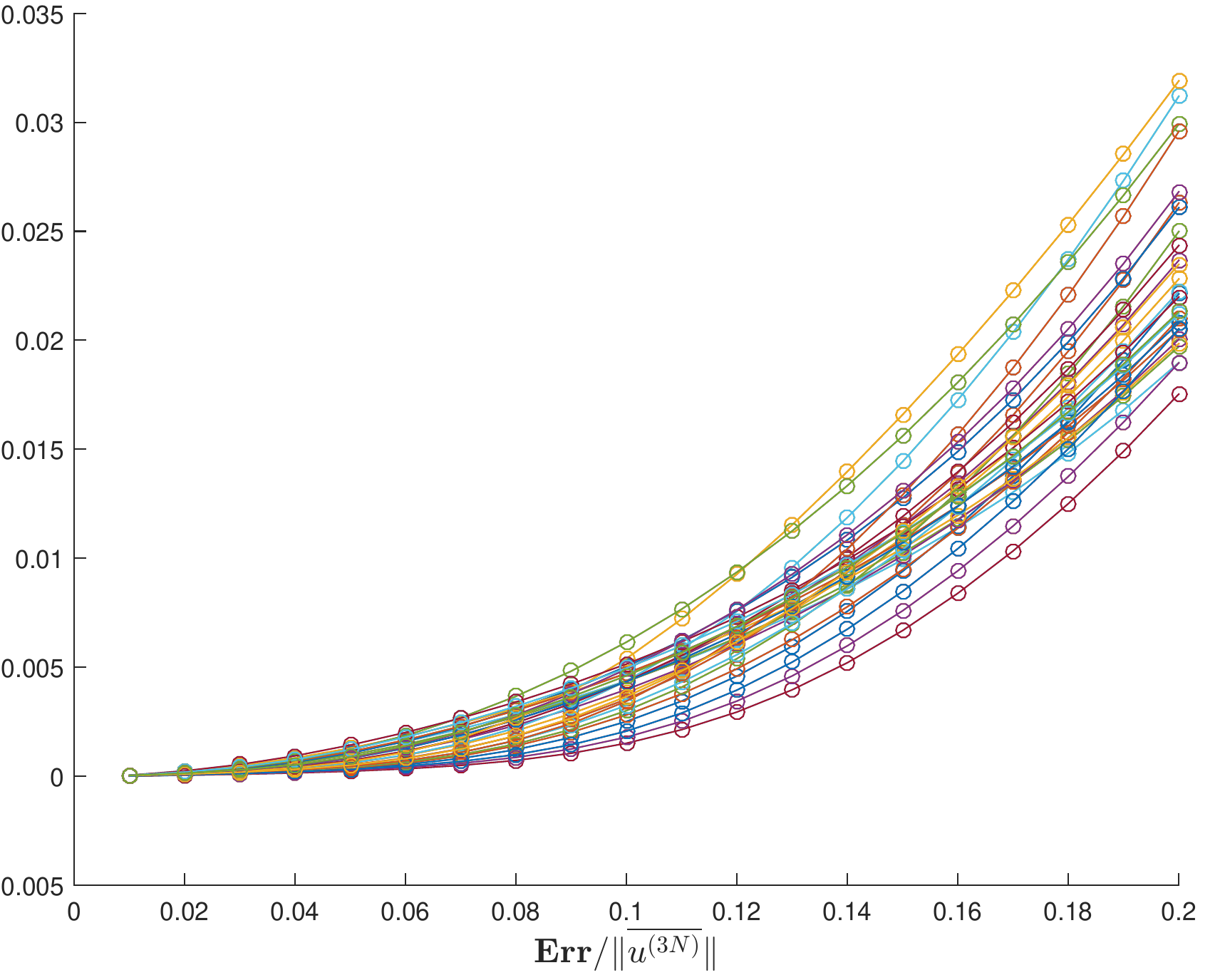}
\caption{Test case 1. Examples of plots of $C(i_1,i_2,i_3)$
against $\mathrm{{\bf Err}}/\|u^{(3N)}\|$.
 } \label{someC}
\end{center}
\end{figure}

\begin{figure}[ht] 
\begin{center}
\includegraphics[scale=.3]{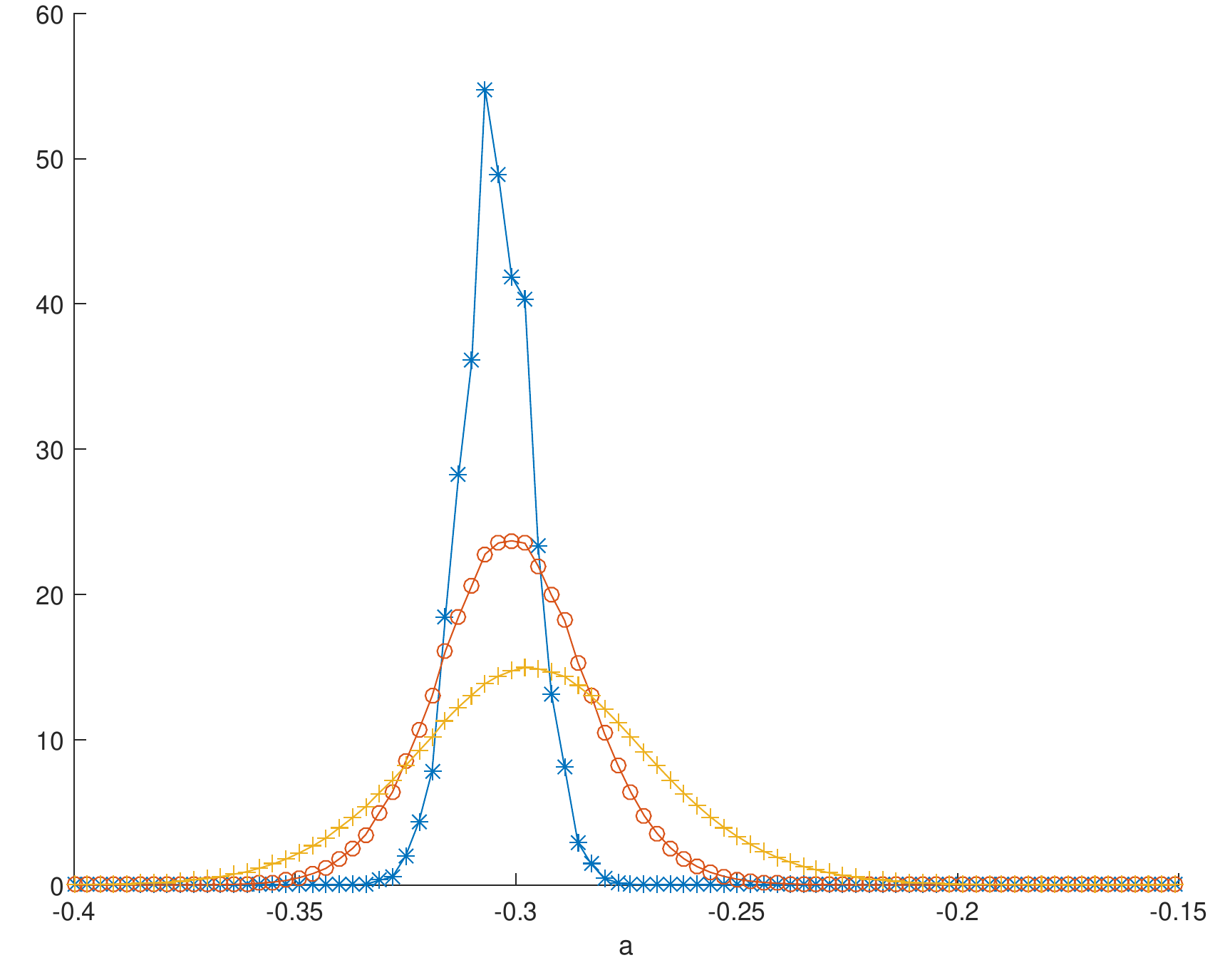}
\includegraphics[scale=.3]{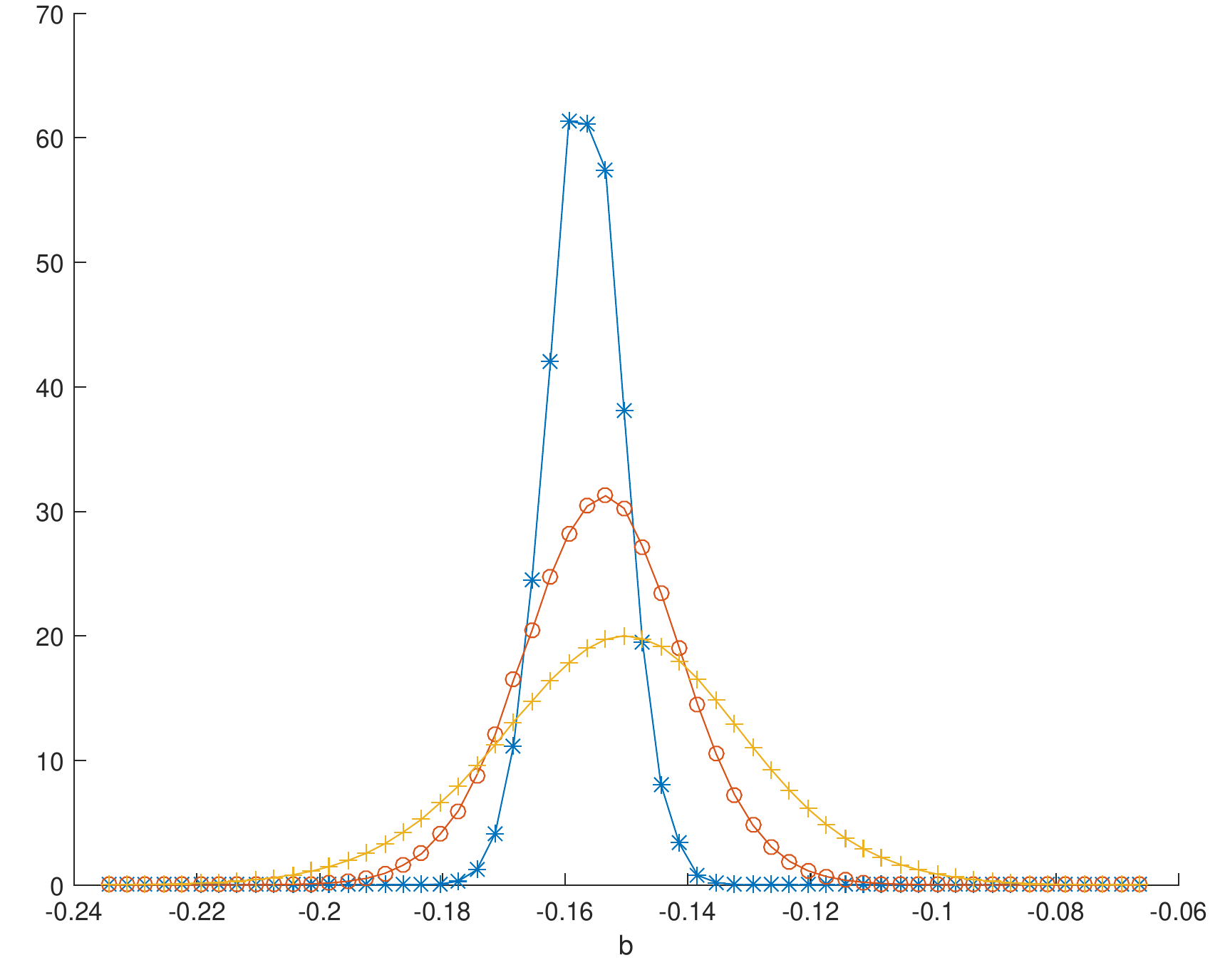}
\includegraphics[scale=.3]{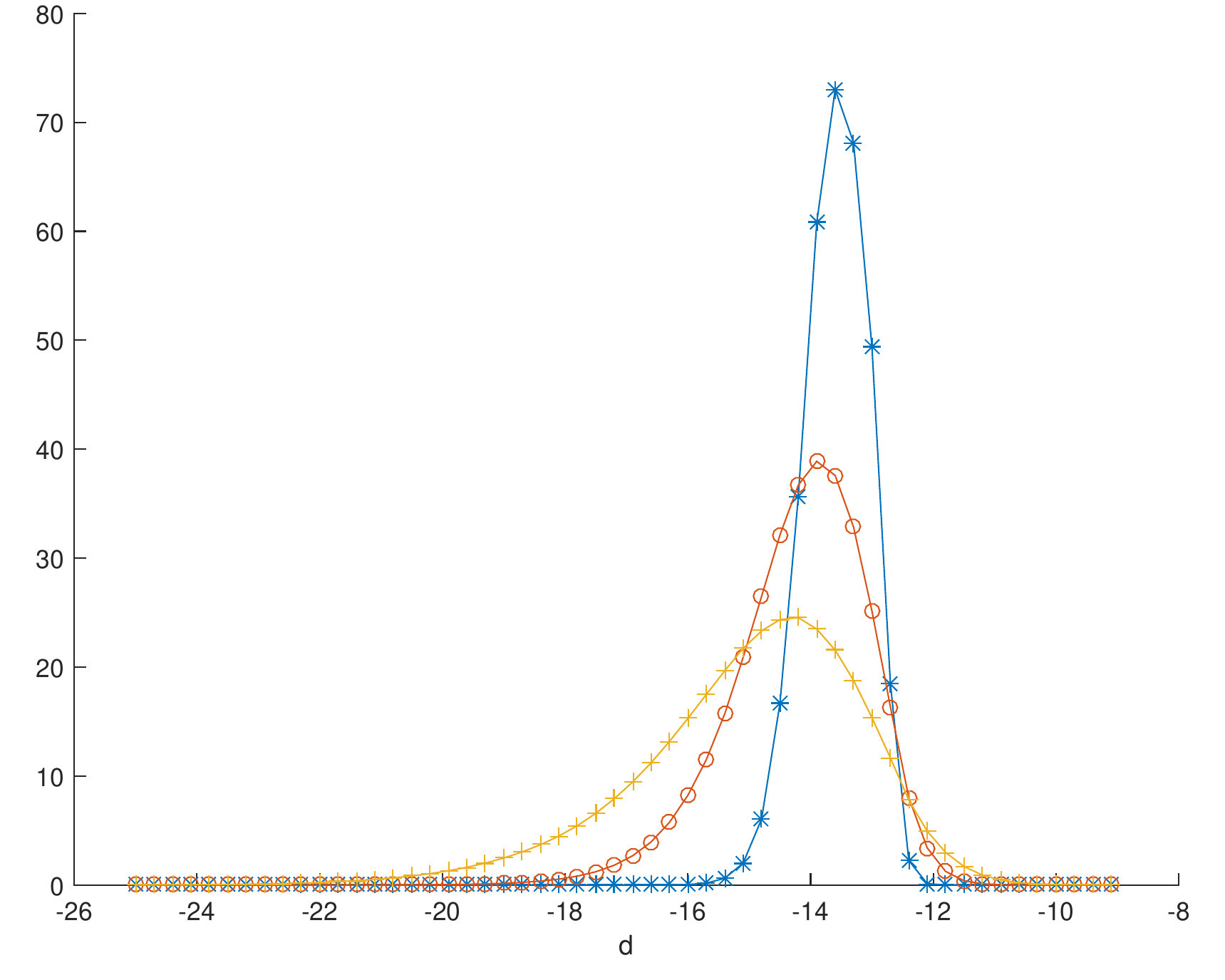}
\caption{Test case 1. Computed marginal distributions for the geometry 
parameters $a$, $b$, and $d$. The blue star curve corresponds to the assumption  $\sigma_{hor}=1, \sigma_{ver}=3$, 
the red circle  curve corresponds to the assumption  $\sigma_{hor}=2, \sigma_{ver}=6$, 
and the orange cross curve corresponds to the assumption  $\sigma_{hor}=3, \sigma_{ver}=9$.
 } \label{dist15}
\end{center}
\end{figure}

In our first example, $\tilde{m}$ is such that $a=-0.3,  \, b=-0.15, \, d=-14 $.
A sketch of the fault $\Gamma$, of the slip field $\tilde{\bh}$,
and the resulting surface measurements $\tilde{\bu}(P_j)$ is shown in 
Figure \ref{geo15}. After surface displacements were computed 
following formula (\ref{new int}), Gaussian noise with 
zero mean was added.  We picked a covariance matrix with
diagonal terms equal to  $(.5 \mbox{ mm})^2 $ for 
horizontal displacements and $(1.5 \mbox{ mm})^2$ for 
vertical displacements.
% We also modeled cross correlation between measurement points by filling 
%out off diagonal terms in the covariance
% matrix with terms that decay exponentially with respect to the distance between 
%measurement points.
% comment on C-reg here 1e-3
In Figure  \ref{someC} we show computed selected values of  $C(i_1,i_2,i_3)$
 near ${\bf C}$
for different values of $\mathrm{{\bf Err}}$, see definitions \ref{defC} and
(\ref{choice C}). 
It is not possible to point to a single preferred value for $\mathrm{{\bf Err}}$, but we
should expect it to be 
 at least twice the standard deviation  of the measurements. 
%
%norm of the diagonal of the covariance matrix.
 %In this present case the relative error $\mathrm{{\bf Err}}/\|\ov{u^{(3N)}}\|$ is shown
%in Figure  \ref{someC} and it should be at least $.03$.
In Figure \ref{someC} we show selected values of $C(i_1,i_2,i_3)$ for 
the relative error $\mathrm{{\bf Err}}/\|u^{(3N)}\|$ between
0.01 and 0.2. 
Since there is no preferred value of $\mathrm{{\bf Err}}$, we choose several possible  
${\bf C}$: in this particular case $k \cdot 10^{-4}$, for $k=1..10$.
In Figure \ref{dist15} we show  
computed marginal distributions for the geometry parameters $a$, $b$, and $d$ for the value
${\bf C} = 10^{-3}$
and three
different assumed values of  
$\sigma_{hor}$ and $\sigma_{ver}$, the standard deviation for horizontal and vertical measurements. 

\subsection{Second test case}
\begin{figure}[ht]
\begin{center}
\includegraphics[scale=.4]{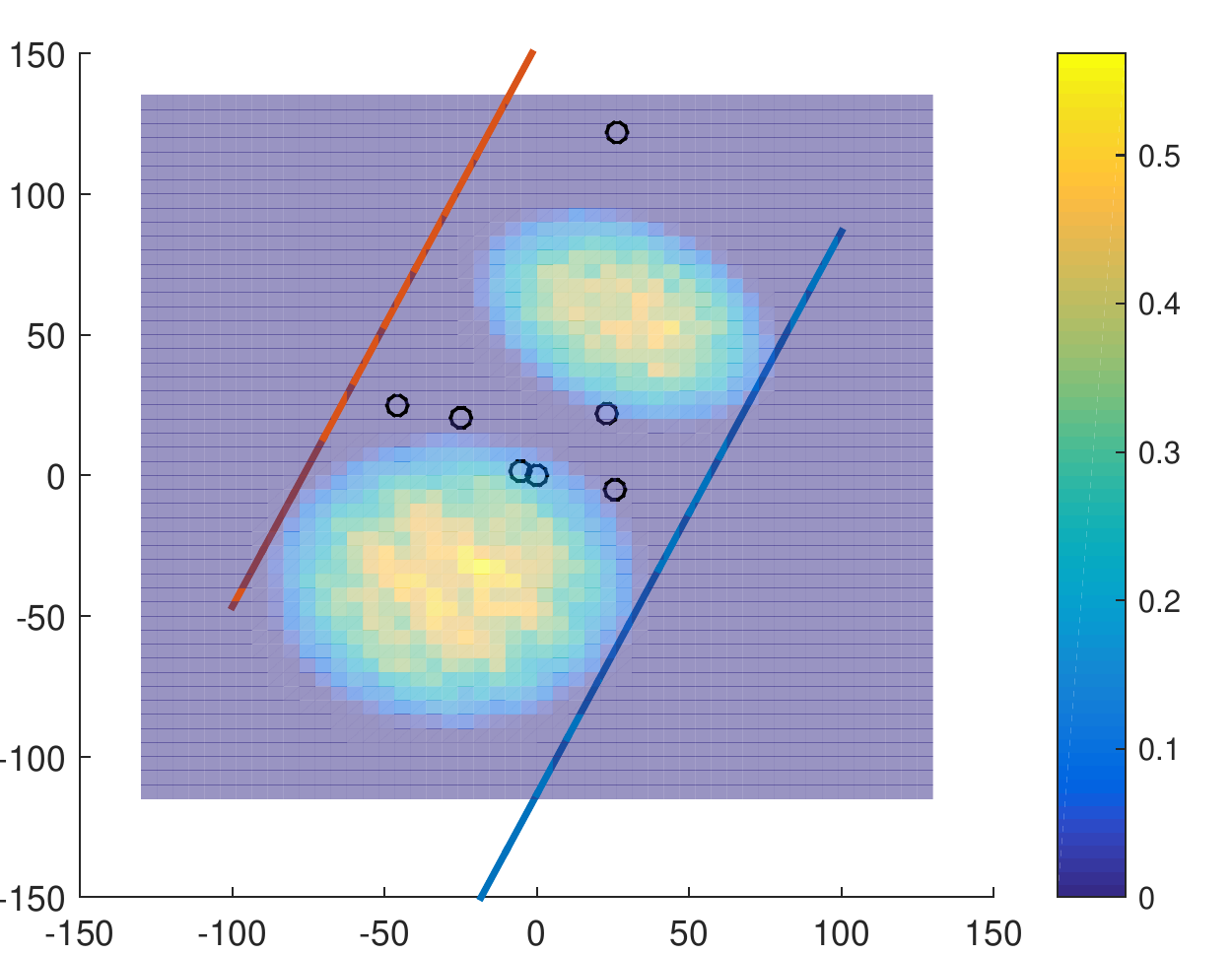}
\includegraphics[scale=.4]{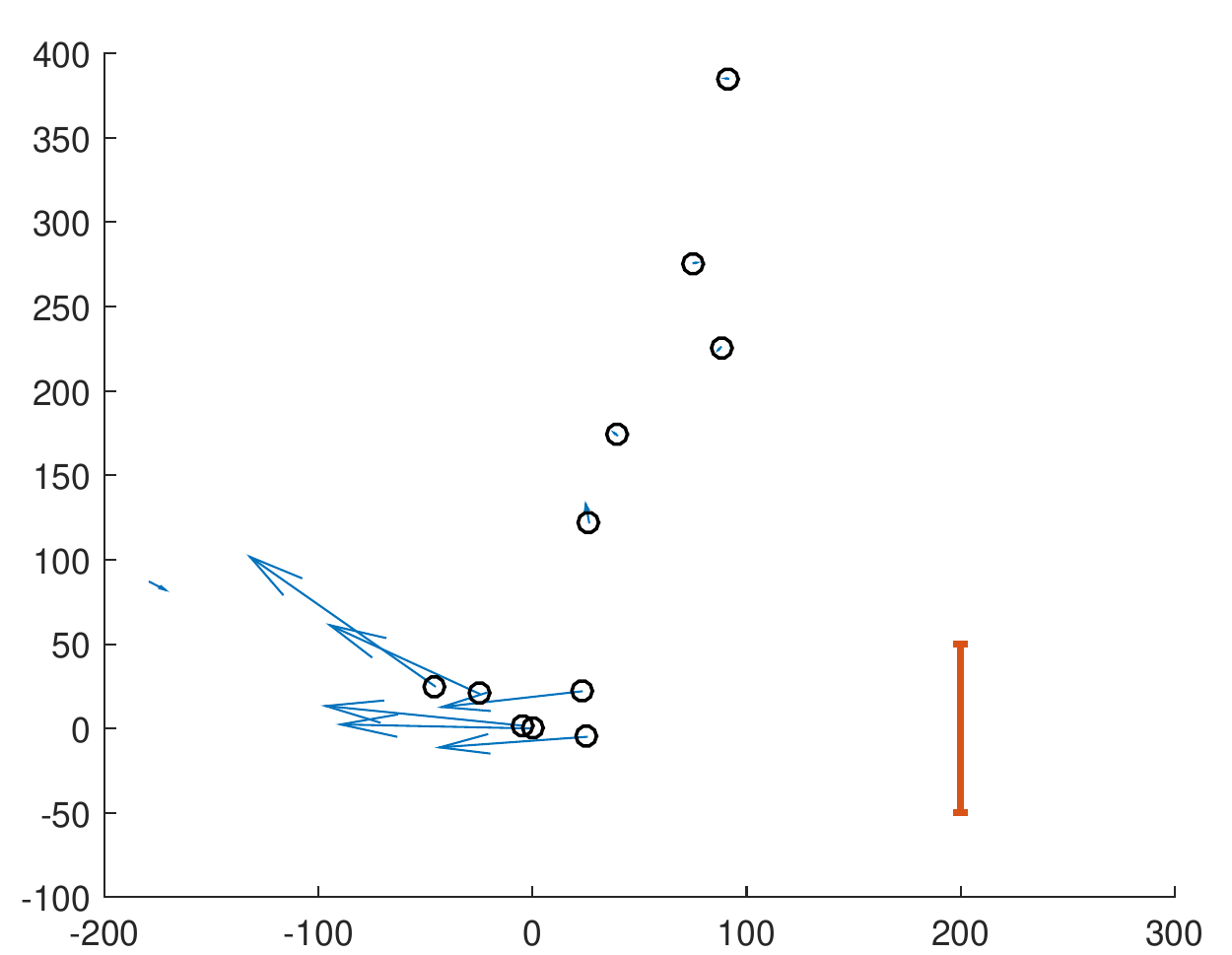}
\includegraphics[scale=.4]{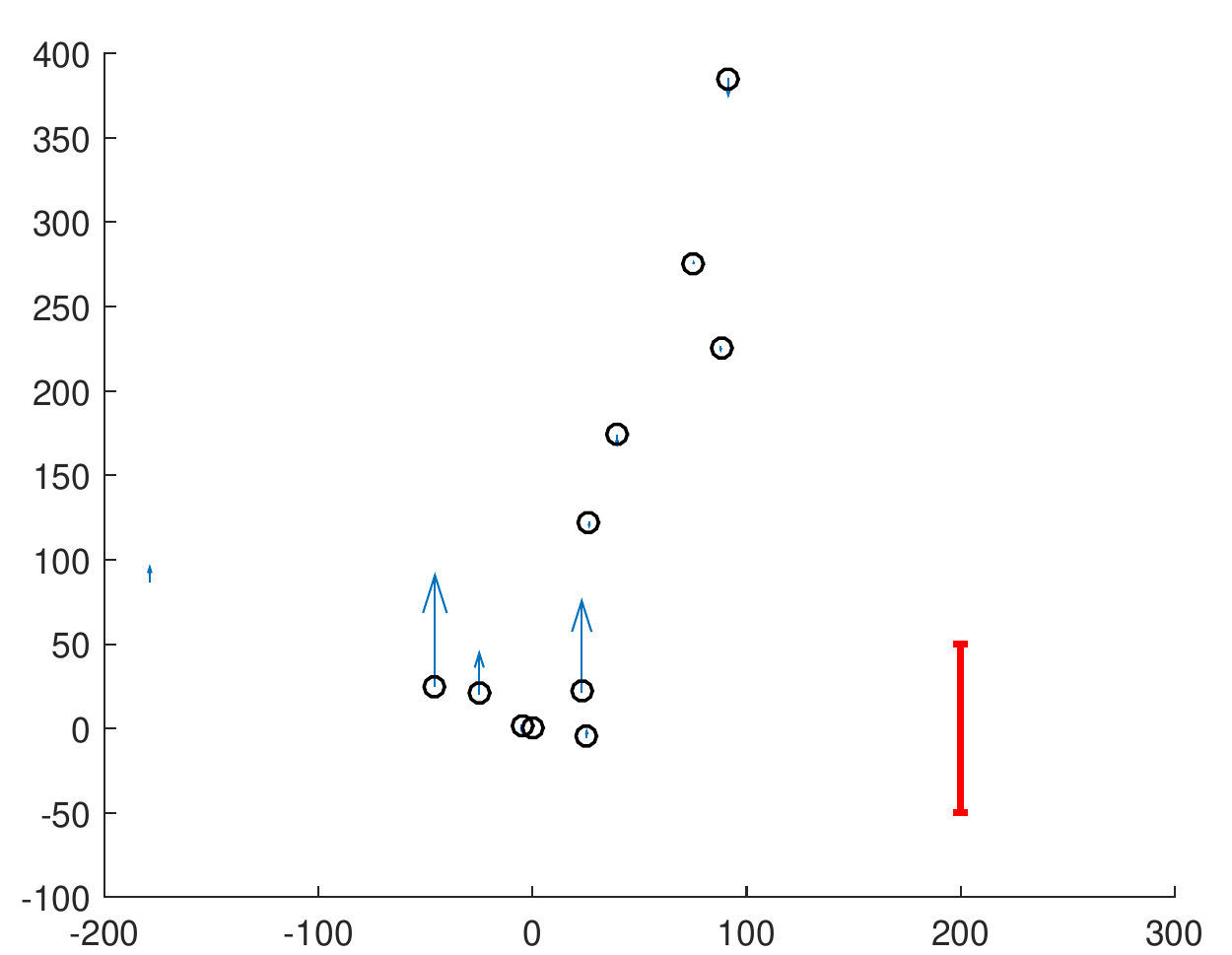}
\caption{Test case 2. Legend is the same as in Figure 
\ref{geo15}.} 
%Top right: reconstruction. Bottom two graphs show horizontal and %vertical displacements,
%data $\bu(P_j)$ in red and reconstructed surface field
%$\tilde{\bu}(P_j)$ in blue  } \label{test2fig}
\label{geomtry2blobs}
\end{center}
\end{figure}

% the so called failure folder
In our second test, $\tilde{m}$ is such that
$a=-0.3,  b=0.15, d=-25 $. 
The slip field for producing the surface data is sketched in Figure \ref{geomtry2blobs}.
This is a more challenging case since this field is non-convex.
In addition, for this combination of geometry and slip field only a few points $P_j$ 
contribute valuable information for the surface displacement field.
In theory, with continuous data on an open set of the surface $x_3=0$ this should not be a problem,
but in practice, with a limited number of observation points our algorithm does not perform as well
as previously.
%
%The secondary zone of slip is obscured by the first one.
%Method is inadequate given how small $N$ is.
The most likely recovered values for $m$ are about
$-0.2, \, 0.1, \, -27$, this is not as close to the correct values as in the previous case.
In Figure \ref{slip 2blobs} we show the reconstructed slip field for this most likely geometry.
Note how one of the two connected components of $\tilde{\bh}$ is better reconstructed than the other one.

\begin{figure}[ht] 
\begin{center}
\includegraphics[scale=.3]{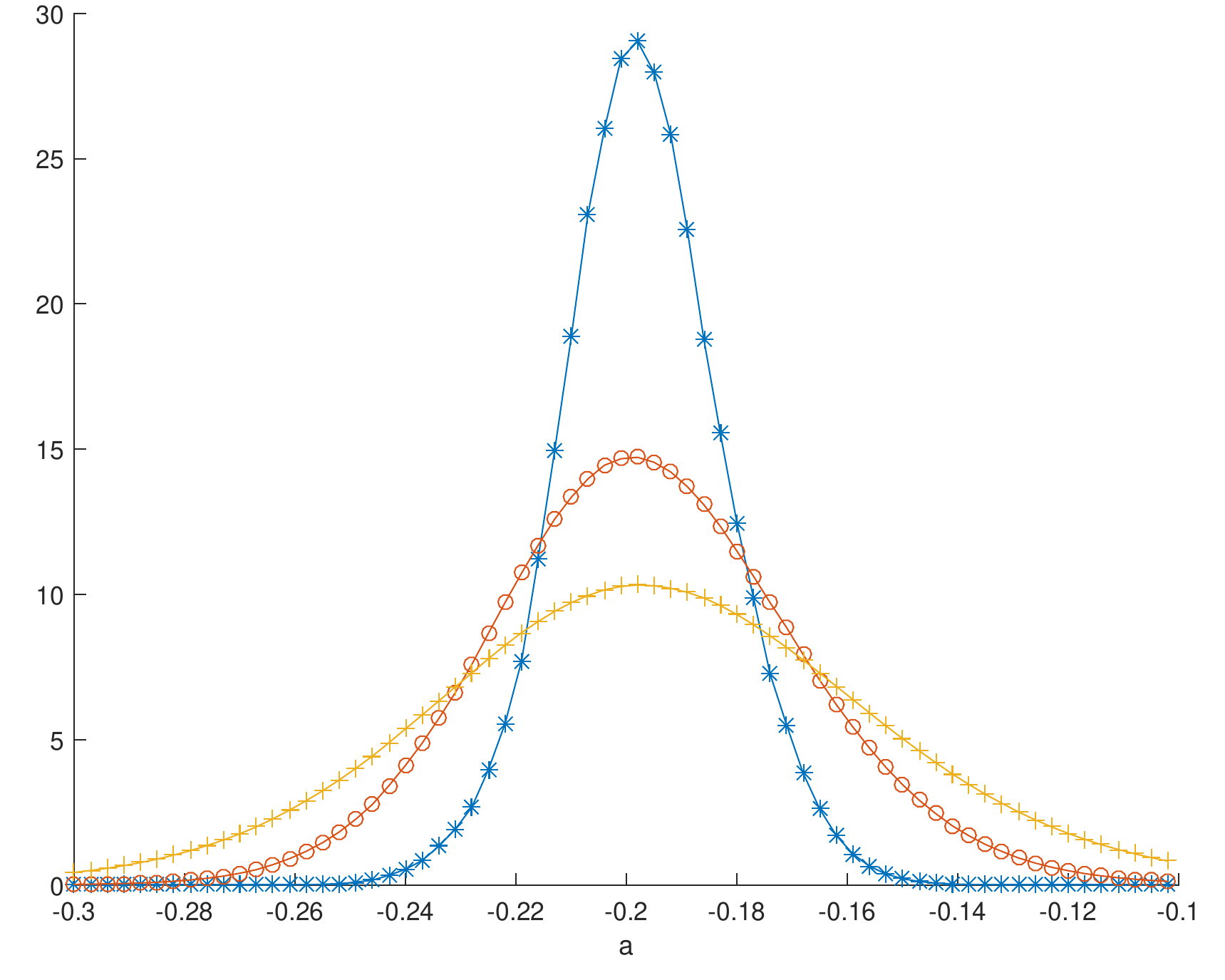}
\includegraphics[scale=.3]{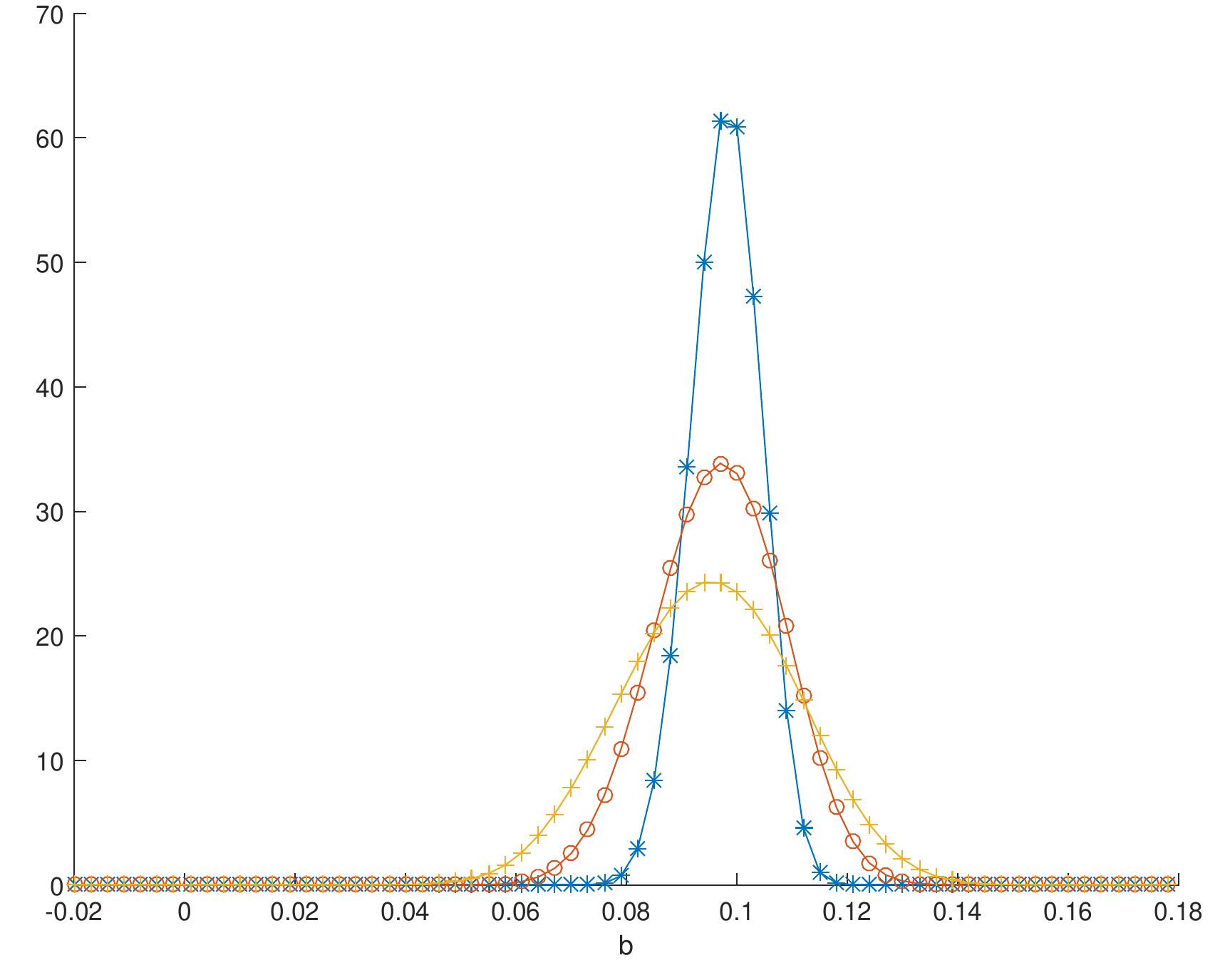}
\includegraphics[scale=.3]{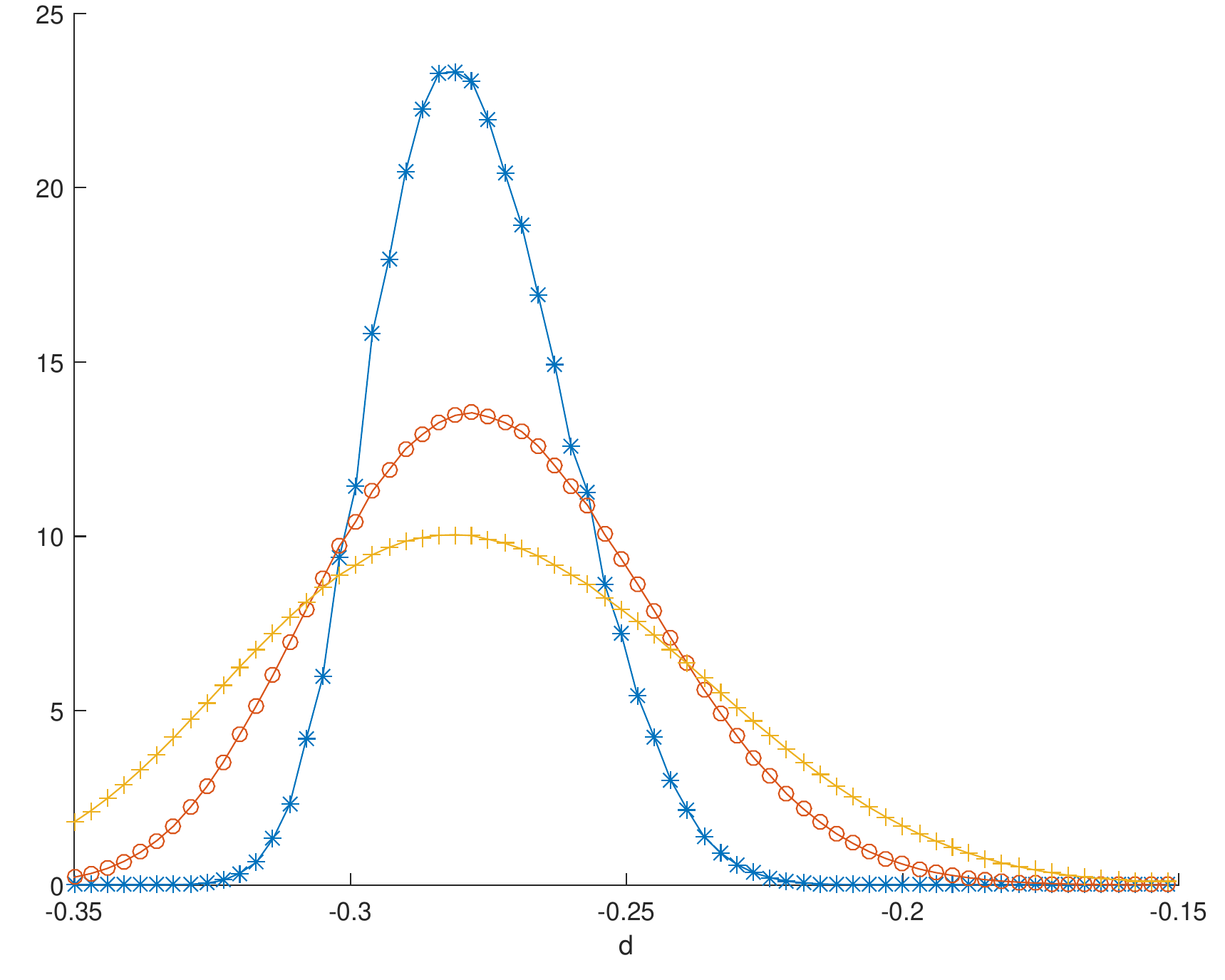}
\caption{Test case 2: computed marginal distributions for the geometry parameters $a$, $b$, and $d$.
Same caption as in Figure \ref{dist15}. 
 } \label{dist2blobs}
\end{center}
\end{figure}

\begin{figure}[ht] 
\begin{center}
\includegraphics[scale=.3]{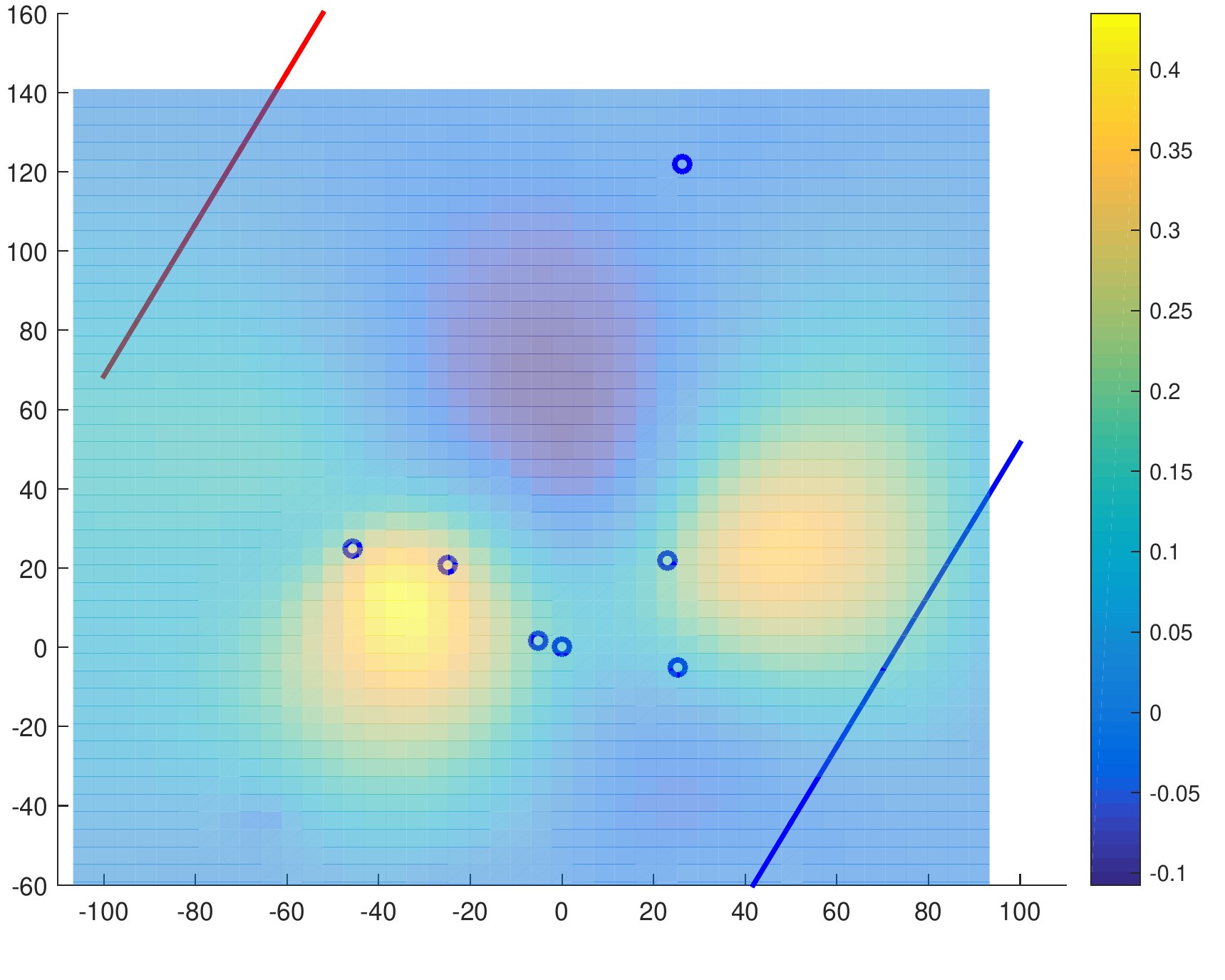}
\caption{Test case 2. Reconstructed slip field.} \label{slip 2blobs}
% done for const_reg=6.e-4;
\end{center}
\end{figure}

\subsection{Third test case}
\begin{figure}[ht] 
\begin{center}
\includegraphics[scale=.4]{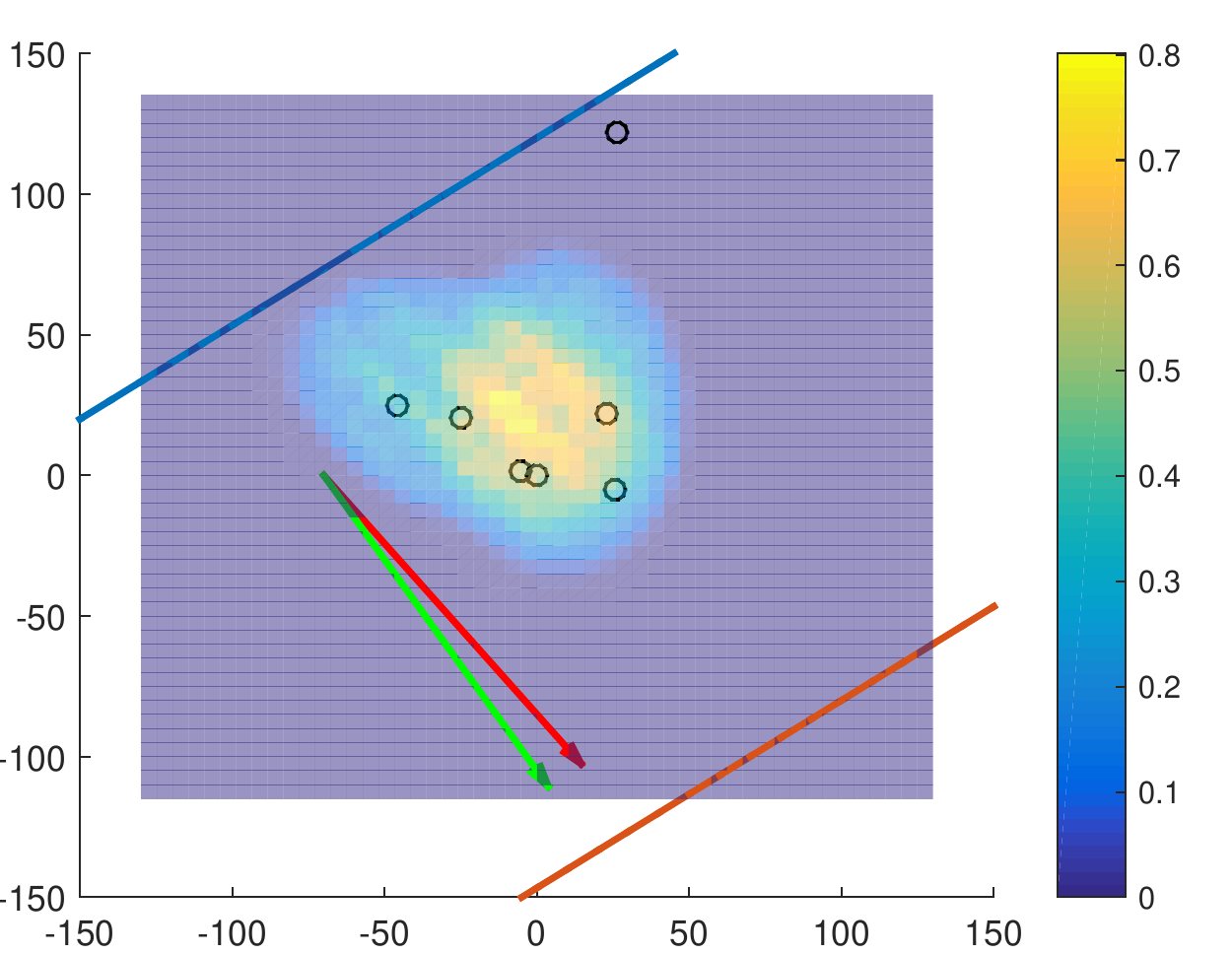}
\caption{Test case 3. Legend is the same as in the top left panel of Figure \ref{geo15}.
The direction of steepest ascent is indicated by the green arrow while the red arrow indicates the
direction of slip.
} \label{geomtryrake}
%Top right: reconstruction. Bottom two graphs show horizontal and %vertical displacements,
%data $\bu(P_j)$ in red and reconstructed surface field
%$\tilde{\bu}(P_j)$ in blue  } \label{test2fig}
\end{center}
\end{figure}

\begin{figure}[ht] 
\begin{center}
\includegraphics[scale=.3]{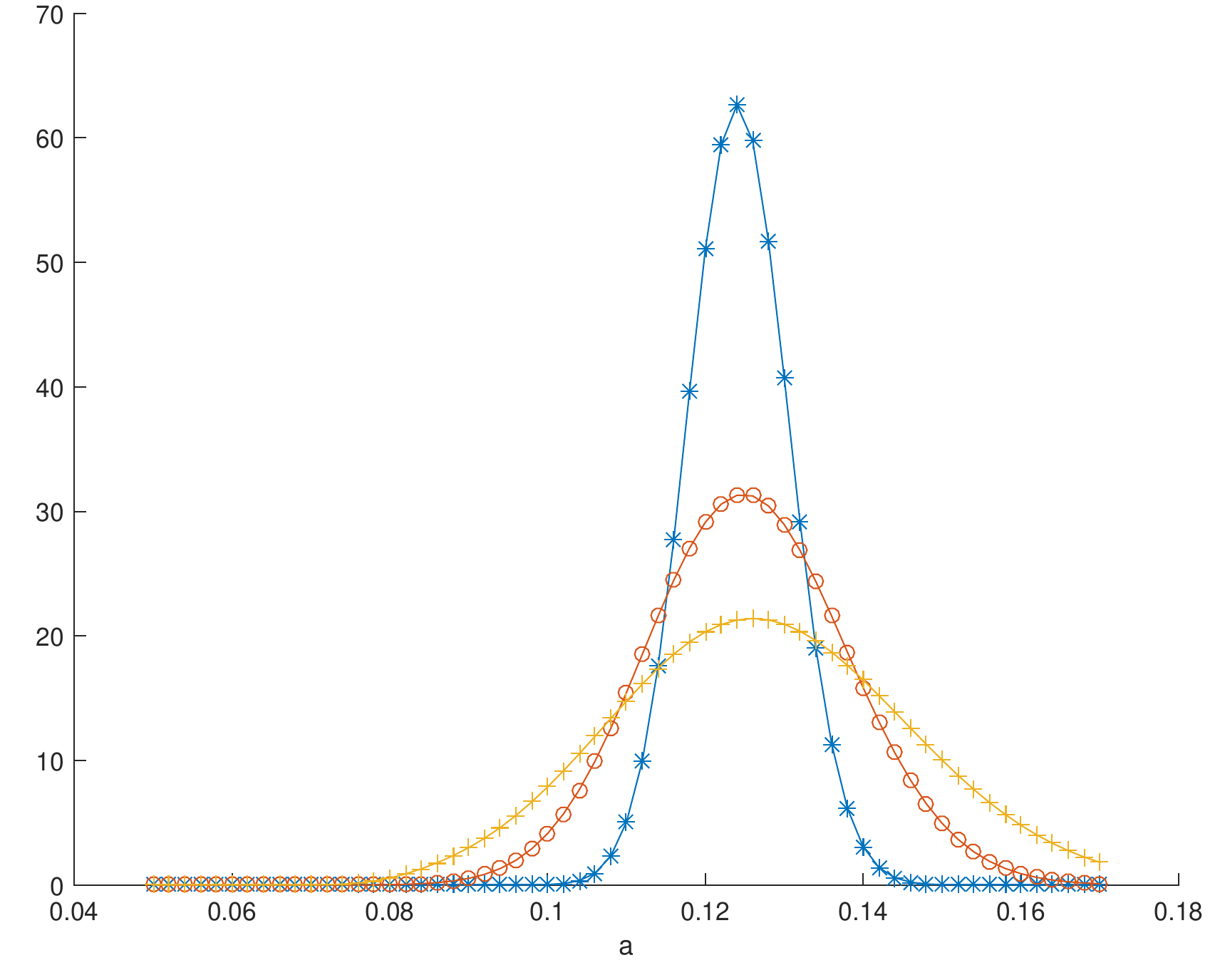}
\includegraphics[scale=.3]{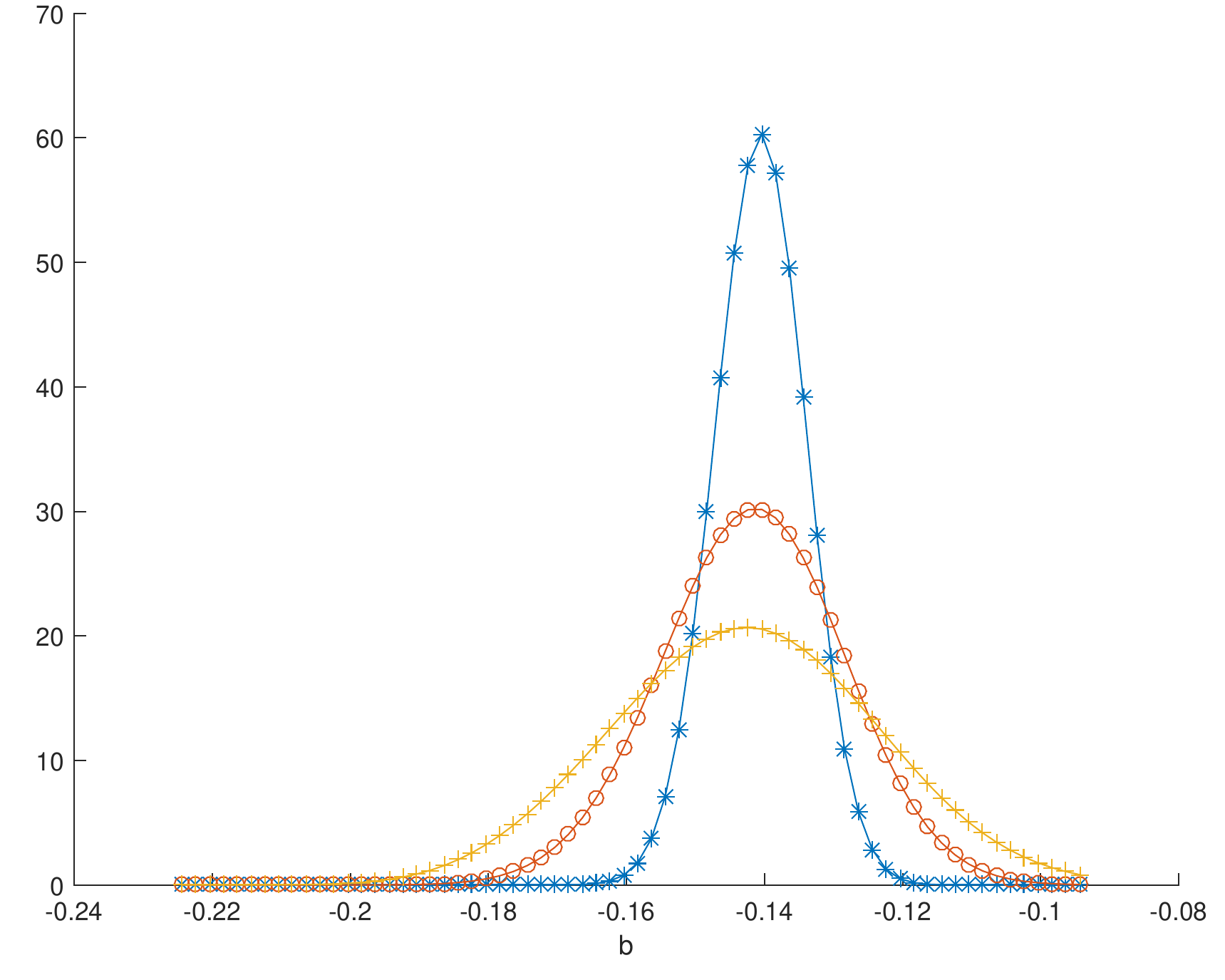}
\includegraphics[scale=.3]{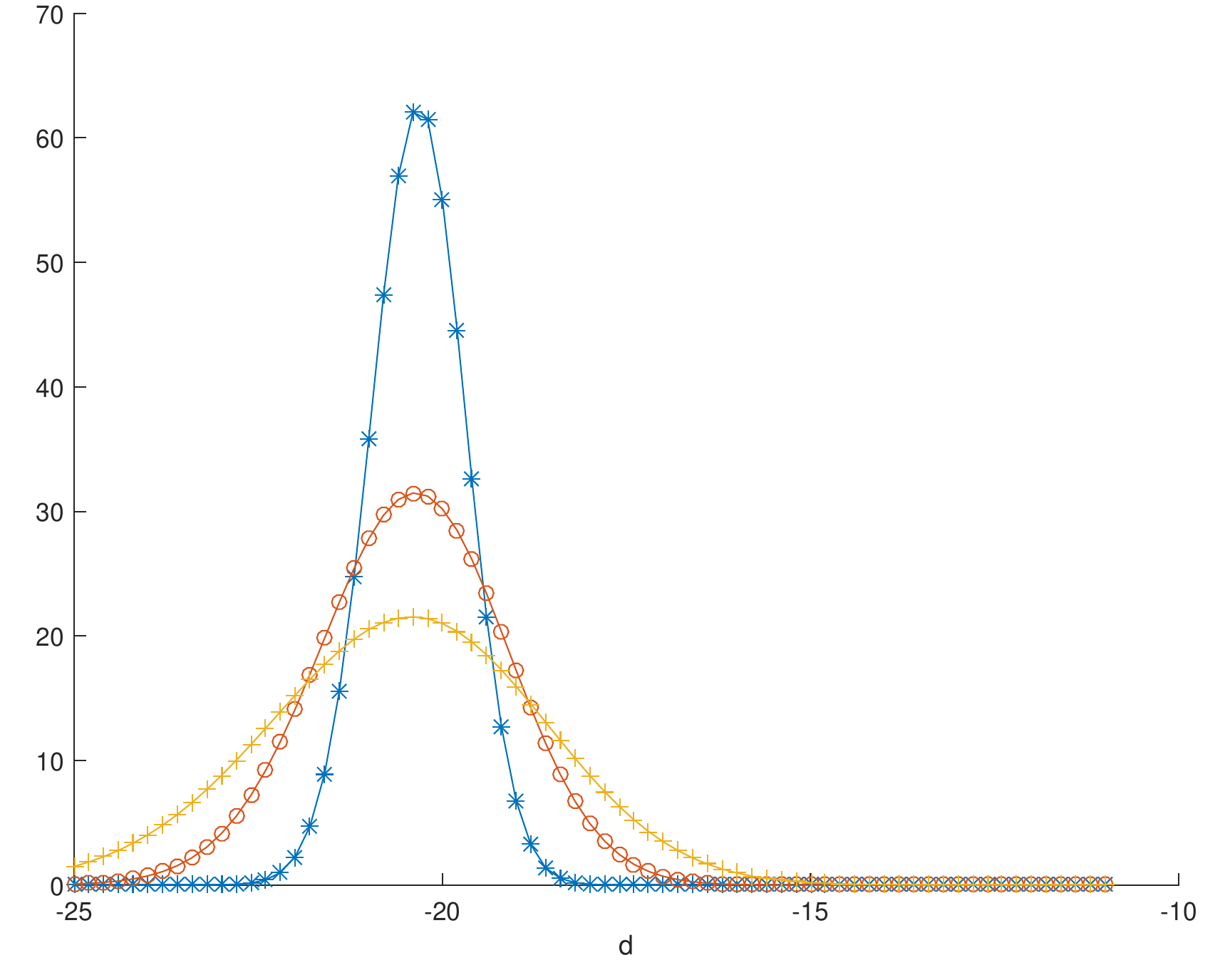}
\caption{Test case 3: computed marginal distributions for the geometry parameters $a$, $b$, and $d$.
Same caption as in Figure \ref{dist15}. } \label{dist975}
\end{center}
\end{figure}

In our third test,  $\tilde{m}$ is such that
$a=0.1,  b=-0.15, d=-24 $. In this case we illustrate how (modest) modeling errors  
may impact the reconstruction algorithm. Here, the direction of slip is {\sl not} in line with 
the direction of steepest ascent, while in the reconstruction step we {\sl wrongly} assume that
these two directions are the same.
These two directions and the fault are sketched in Figure \ref{geomtryrake}.
In addition, noise was added to the surface measurements as in the previous two cases.
In Figure \ref{dist975} we show computed marginal distributions
 for the geometry parameters $a$, $b$, and $d$.
The computed maximum likelihood for $m$
 are achieved at .12, -.14, -20, so in this "wrong model" case these values are not as close to
the original values that were used to produce data as they were in the first test case.
%that s 9point75


\subsection{Application to the case of measured surface displacements during the 2007 SSE
in Guerrero, Mexico}
\begin{figure}[ht] 
\begin{center}
\includegraphics[scale=.3]{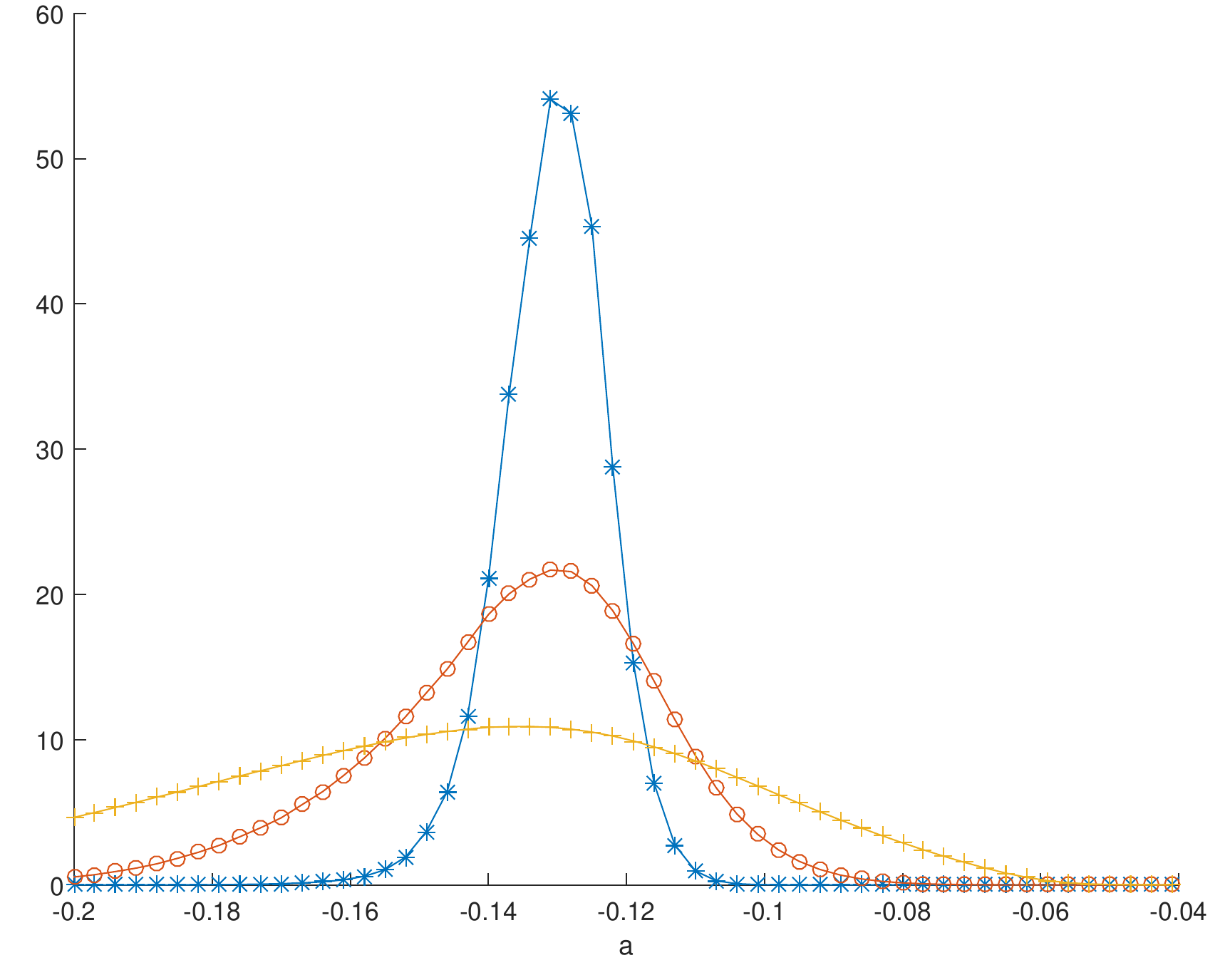}
\includegraphics[scale=.3]{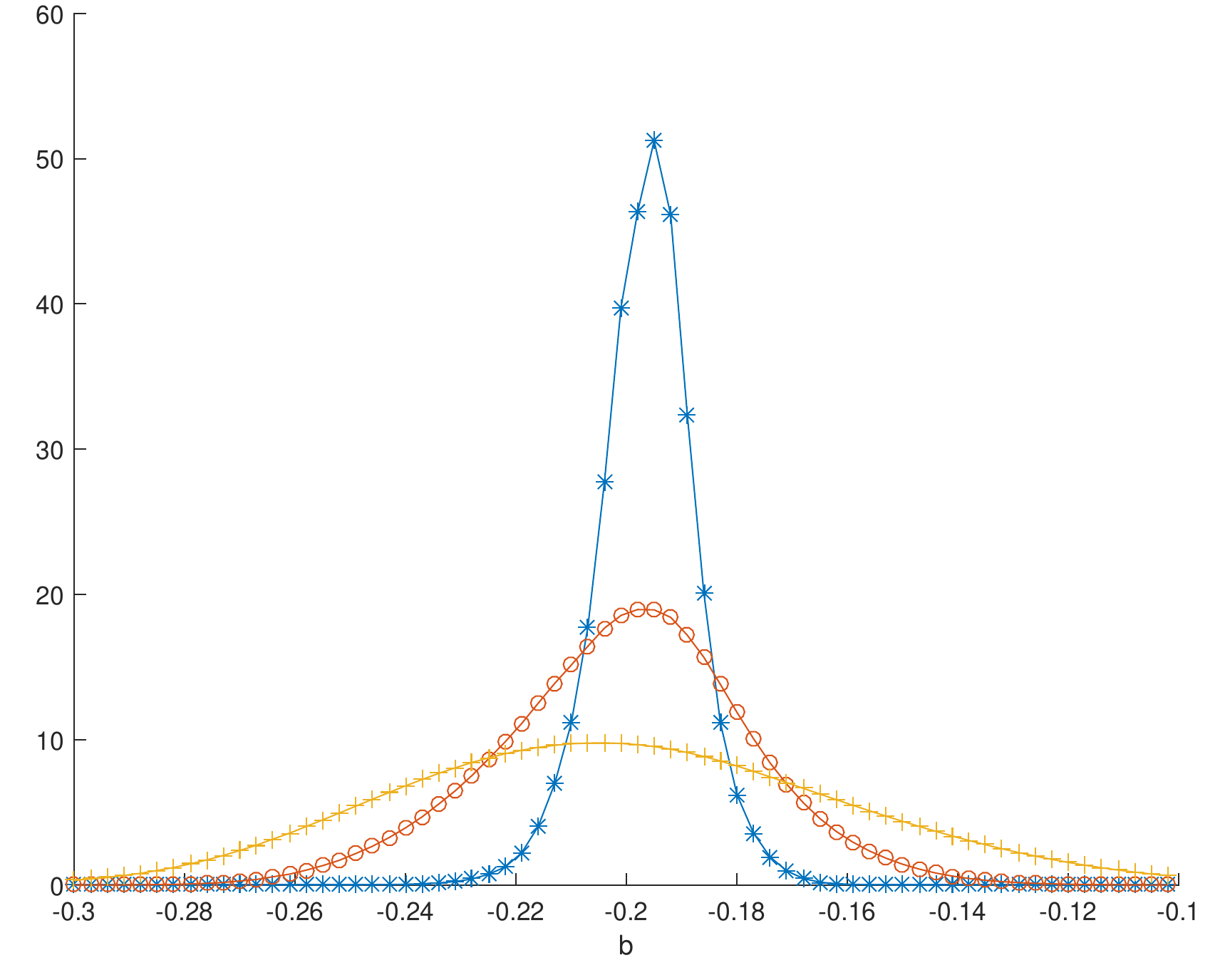}
\includegraphics[scale=.3]{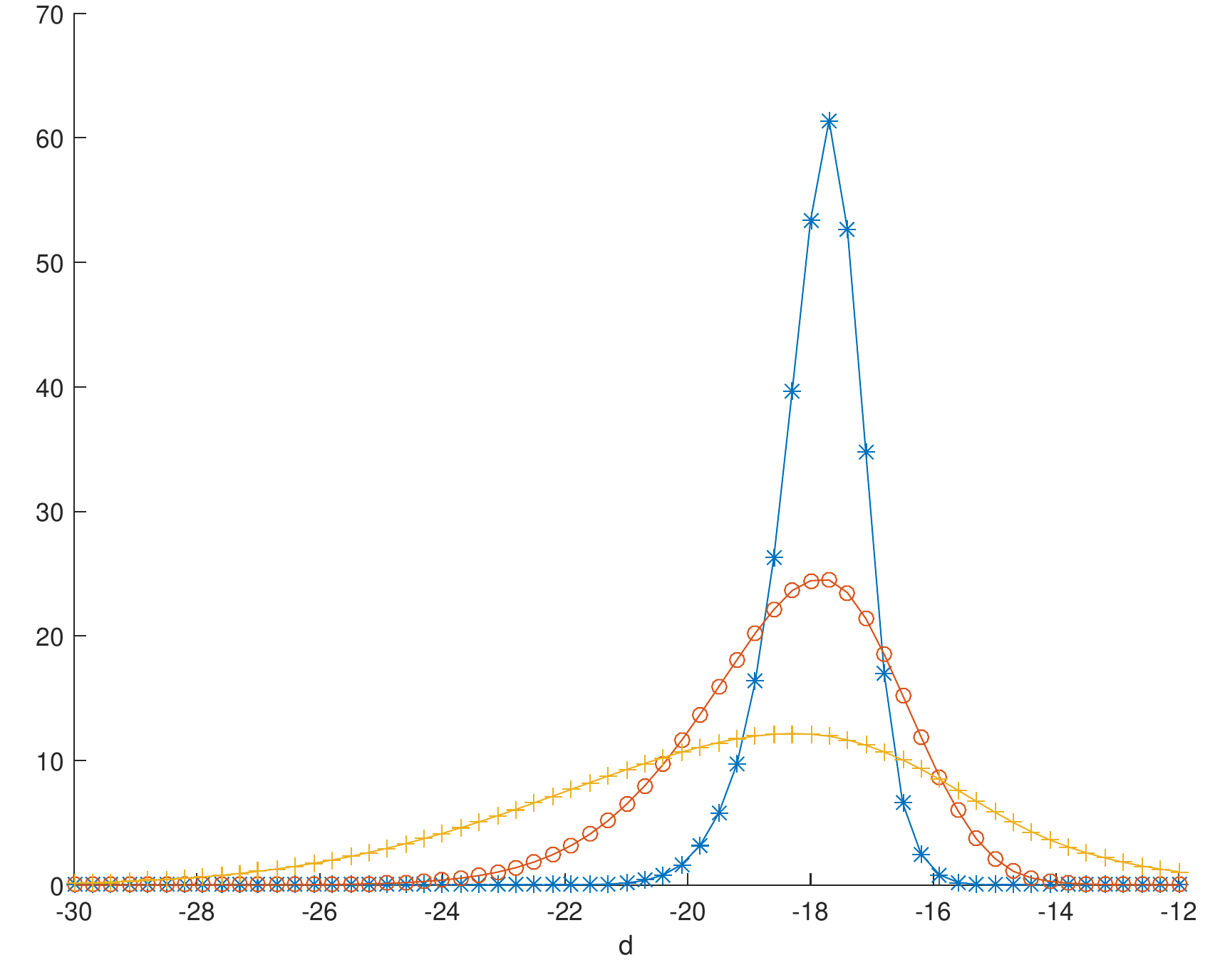}
\caption{The 2007 Guerrero SSE.
Computed marginal distributions for the geometry parameters $a$, $b$, and $d$.
The blue star curve corresponds to the assumption that $\sigma_{hor}=.5, \sigma_{ver}=1.5$, 
the red circle  curve corresponds to the assumption that $\sigma_{hor}=1, \sigma_{ver}=3$, 
and the orange cross curve corresponds to the assumption that $\sigma_{hor}=2, \sigma_{ver}=6$.
  } \label{distmeas}
\end{center}
\end{figure}

\begin{figure}[ht] 
\begin{center}
\includegraphics[scale=.3]{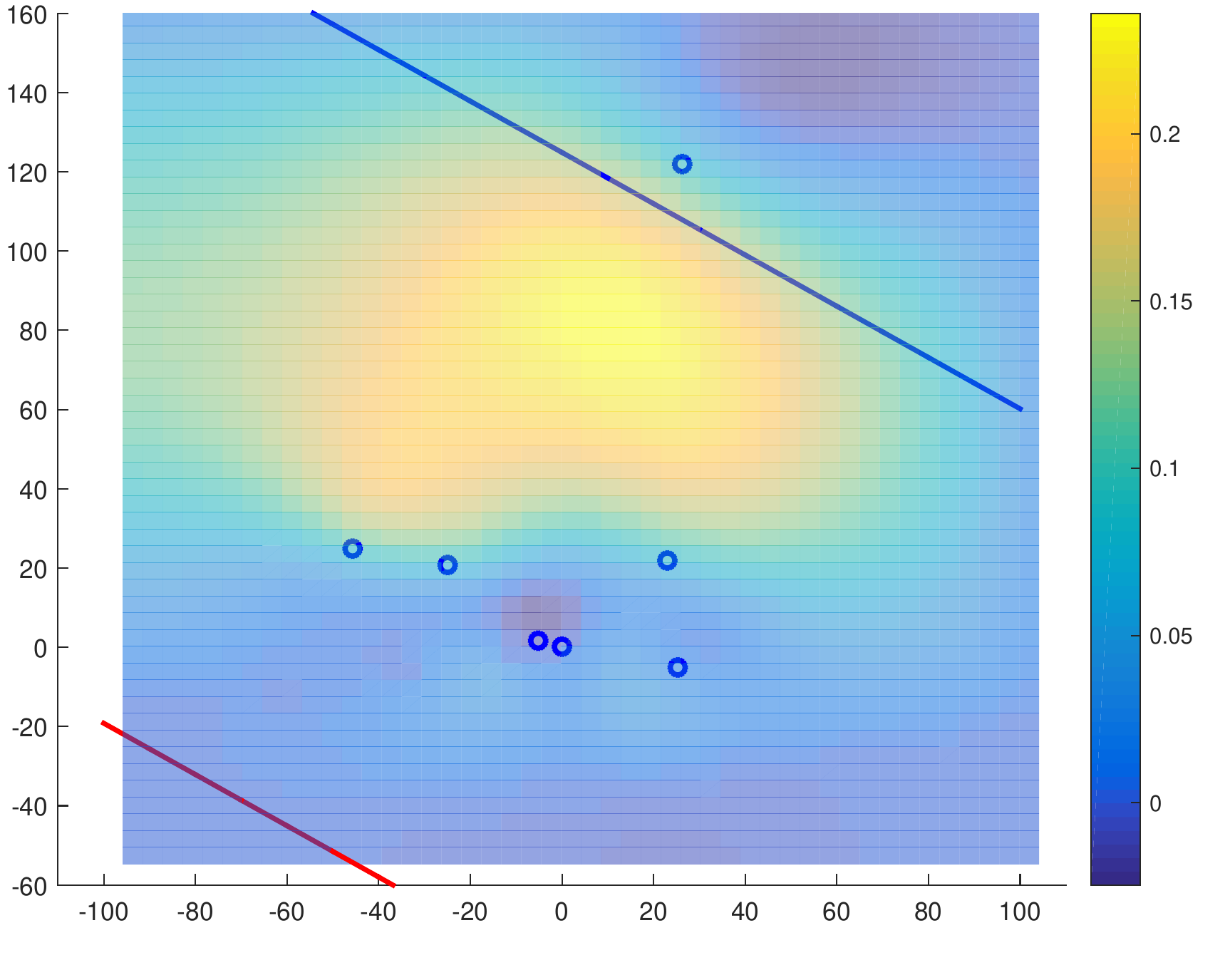}
\includegraphics[scale=.3]{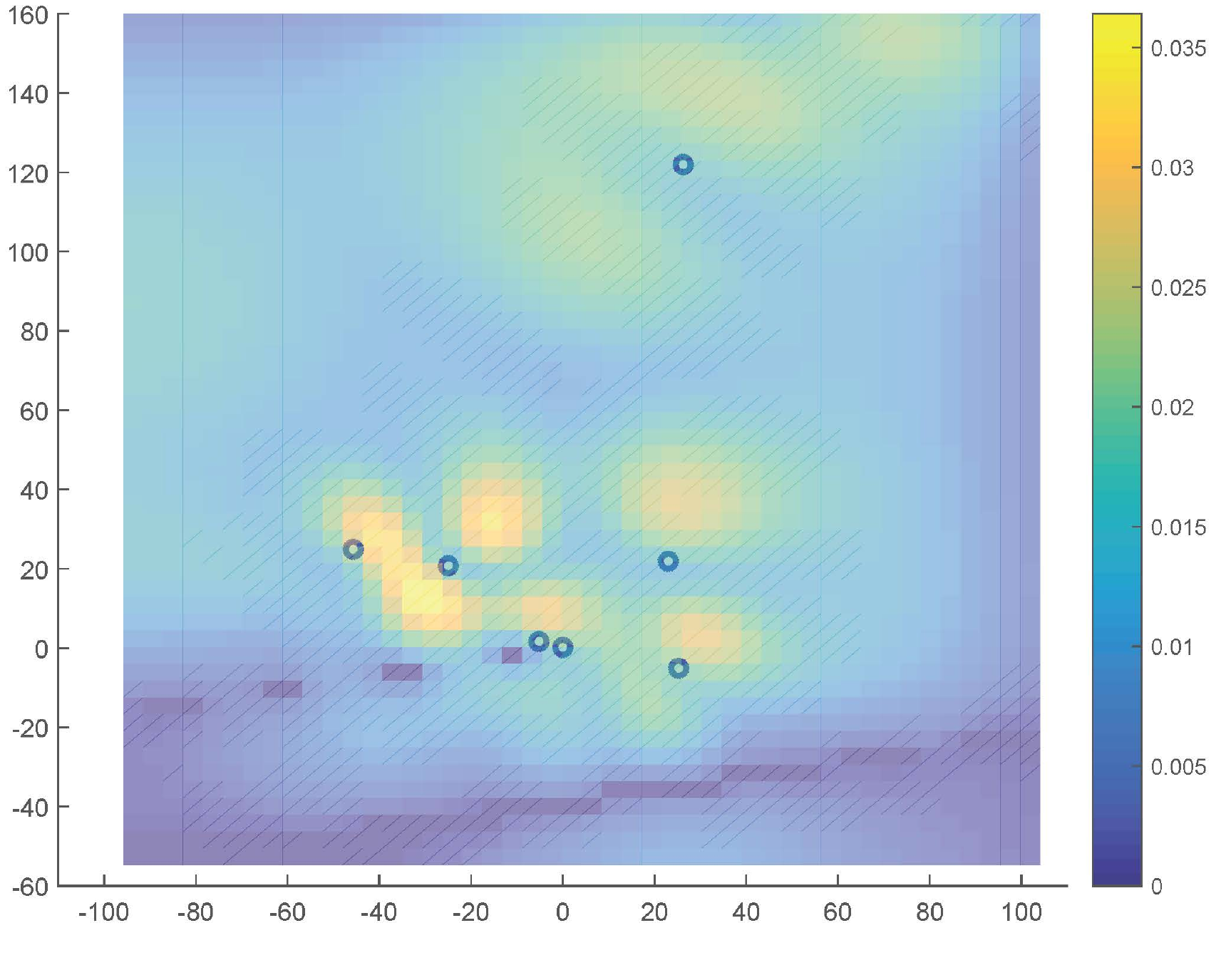}
\caption{Computed average slip (left) and standard deviation (right) for the Guerrero 2007 SSE. Note the change of scale for the color bars
between the two figures.} \label{slip meas}
% done for const_reg=6.e-4;
\end{center}
\end{figure}

% WARNING!!
%  in these calculations the sub-horizontal part was considered
%   the part that is too close to the surface was always set to zero
%   so the inverse of matrix in algorithm had to be - recomputed at each step
% (it depends on geometry)
We now show the most interesting case as far as applications are concerned.
We start from measurements relative to the 2007 SSE
in Guerrero, Mexico, which were processed as described earlier: both $\bu(P_j)$ and 
standard deviation on these measurements were estimated.
We show in Figure \ref{distmeas} 
computed marginal distributions for the geometry parameters $a$, $b$, and $d$
for the constant  ${\bf C}$ set to $6 \cdot 10^{-4}$.
Next we fix (approximate) most likely values for the geometry parameters $a, b, d$ to
$-.13, \, -.19, \, -18$, and we compute expected slip on the fault and standard deviation:
results are shown in Figure \ref{slip meas}.
Here we need to point out that once the geometry of the fault is fixed, we only need
to solve a {\sl linear} stochastic inverse problem: this is rather trivial since there is a linear
relationship between the covariance matrix of the data and the covariance matrix 
of the slip on the field.\\
In the case of measured data, we can only validate our calculation by comparing 
our reconstructed fault to those offered by earlier studies:
see \cite{Kostoglodov, Pacheco, Suarez} for the geometry of the fault (these studies were based
on seismicity and gravity),
\cite{Radiguet, Radiguet2} for the profile of the slip on the fault, and
\cite{volkov2017determining, volkov2017reconstruction} for combined (deterministic) studies of 
simultaneous reconstruction of geometry and slip fields.
In Figure \ref{slip meas},  the computed line with depth $x_3=-2$ on the fault is shown in red.
Note how close to the Middle American Trench sketched in Figure \ref{geography} this line is.
With $\sigma_{hor}=1, \sigma_{ver}=3$, the standard deviations for $m$ 
are 
$0.020, \, 0.023, \, 1.7$, so the depth below ACAP is approximately between 16 and 20 km, and
50 km in the direction of steepest descent the depth is between 27 km and 33 km 
(these are plus or minus 1 standard deviation intervals).
This is comparable and rather on the high side of depths found in other studies, see Figure 10 in
\cite{volkov2017reconstruction}.

%d=-.18071;
%b=-.19158;
%a=-.12398;

\clearpage 
\section{Appendix}
The following two lemmas are needed for the proof  of Theorem \ref{fdisc con}.
\begin{lem} \label{bhdisc bounded}
Let $\bh^{disc}_{m,C}$ be the unique point where 
$F^{disc}_{m,C}$ achieves its minimum
	in ${\cal F}_p$.
There are two positive constants $C_0$ and $C_1$ such that 
$\bh^{disc}_{m,C}$ is uniformly bounded in $H^1_0(R)$ for all $m$ in $B$, $p$ in $\NN$, $N$ in $\NN$ and $C$ such that 
$C_0 N^{-\beta} < C <C_1$.
\end{lem}
\textbf{Proof:} \\
With $\tilde{\bh}$  as in the statement of Theorem \ref{cv th},
let $\tilde{\bh}_p$ be the orthogonal projection of $\tilde{\bh}$
on ${\cal F}_p$.
As
\bea
F^{disc}_{m,C} (\bh_{m,C}^{disc}) \leq F^{disc}_{m,C} (\tilde{\bh}_p),
\eea
we have 
\bean \label{bounded1}
C \int_R |\bh_{m,C}^{disc}|^2 \leq \sum_{j=1}^N  C'(j,N)  | {\cal C}^{-\f12}(A_{m} 
\tilde{\bh}_p  
-  \tilde{\bu} )(P_j)|^2 + C \int_R |\tilde{\bh}_p|^2.
\eean
Since $\| \tilde{\bh}_p\| \leq \| \tilde{\bh}\|$,
given that $A_m $ is continuous in $m$, $B$ is compact,
 and $A_m$ continuously maps $H^1_0(R)$ into smooth functions on $V$,
 by (\ref{quad rule}),
\bea
\sum_{j=1}^N  C'(j,N)  | {\cal C}^{-\f12}(A_{m} \tilde{\bh}_p  
-  \tilde{\bu} )(P_j)|^2 = O(N^{-\beta}),
\eea
thus
\bean \label{bounded2}
\int_R |\bh_{m,C}^{disc}|^2 = O(1 + C^{-1} N^{- \beta}),
\eean
uniformly   for all $m$ in $B$, $p$ in $\NN$, $N$ in $\NN$ and $C>0$.

%\begin{lem} \label{bhdisc converges}
%Fix $\epsilon >0$. 
%Let $\bh_{m,C}^{disc}$ be as in Lemma \ref{bhdisc bounded}, and $\bh_{m,C}$ be in $H^1_0(R)$ such that 
%$F_{m,C} (\bh_{m,C} )$ minimizes $F_{m,C}$ over $H^1_0(R)$. 
%There is a $p_0$ in $\NN$, $N_0$ in $\NN$, and two positive constants
%$C_0$ and $C_1$ such that for all $p \geq p_0$, $N \geq N_0$ and 
%$C$ such that $C_0 N^{-\beta} < C <C_1$,
%\bean \label{F con}
%|F_{m,C} (\bh_{m,C} ) - F_{m,C}^{disc} (\bh_{m,C}^{disc} ) | \leq \epsilon.
%\eean
%\end{lem}
%\textbf{Proof:} \\
%Since $\bigcup_{p=1}^\infty {\cal F}_p$ is dense in $H^1_0(R)$ and the sequence ${\cal F}_p$
%is increasing for inclusion, there is a $p_0$ in $\NN$ 
%and an $N_0$ in $\NN$ such that for all for all $p \geq p_0$ and  $N \geq N_0$,
%\bea
%F^{disc}_{m,C} (\bh^{disc}_{m,C}) \leq \f{\epsilon}{2} + \min_{H^1_0(R)} F^{disc}_{m,C} 
%\leq \f{\epsilon}{2} +F^{disc}_{m,C} (\bh_{m,C}) \leq \epsilon + F_{m,C} (\bh_{m,C}).
%\eea
%Conversely, 
%thanks to Lemma \ref{bhdisc bounded}, 
%\bea
%F_{m,C} (\bh_{m,C}) \leq F_{m,C} (\bh^{disc}_{m,C}) \leq \epsilon +
%F^{disc}_{m,C} (\bh^{disc}_{m,C}),
%\eea
%for  $N \geq N_0$ and 
%$C$ such that $C_0 N^{-\beta} < C <C_1$, where $C_0$ and $C_1$  are two positive constants.

\begin{lem} \label{almost stable}
Assume that 
	$\tilde{\bu} = A_{\tilde{m}} \tilde{\bh}$ for
	some $\tilde{m}$ in $B$ and some 
	$\tilde{\bh}$ in  $H_0^1(R)$.
	Fix $m$  in $B$ such that $m \neq \tilde{m}$ and $M>0$. Set
	\bea
	 \gamma = \inf_{\bg \in H^1_0(R), \| \bg \| \leq M} \int_V | {\cal C}^{-\f12}(A_m \bg  - \tilde{\bu})|^2.
	\eea
	Then $\gamma >0$.
\end{lem}
\textbf{Proof}:\\
Arguing by contradiction, assume that $\gamma =0$. Then there is a sequence $\bg_n$ in $H^1_0(R)$
such that $\| \bg_n\| \leq M$ and $ \int_V | {\cal C}^{-\f12}(A_m \bg_n  - \tilde{\bu})|^2 $ converges to zero as
$ n \ri \infty$. A subsequence of $\bg_n$ is weakly convergent in $H^1_0(R)$ to some $\bh^*$. 
It will still be  denoted by $\bg_n$  for the sake of simpler notations.
As the operator $A_m$ is compact, we find at the limit that
$ \int_V | {\cal C}^{-\f12}(A_m \bh^*  - \tilde{\bu})|^2  =0$. Since $m \neq \tilde{m}$, this contradicts 
uniqueness  Theorem \ref{uniq1}.\\

\textbf{Proof of Theorem \ref{fdisc con}}:\\
Arguing by contradiction, assume that there exist an $\eta >0$ and  three sequences 
$N_n$ in $\NN$, $p_n$ in $\NN$, and $C_n$ in $(0, 1)$ such that 
$N_n \ri \infty$, $p_n \ri \infty$ and $C_n \ri 0$
while $C_n^{-1} N_n^{- \beta}$ is bounded above and 
denoting 
 $$
f^{disc}_{C_n} (m_n) = \min_{m \in B} \min_{\bg \in {\cal F}_{p_n}}\sum_{j=1}^{N_n} C'(j,N_n) 
| {\cal C}^{-\f12}(A_{m} \bg - \tilde{\bu} ) (P_j) |^2 
+ C_n \int_R |\nabla \bg|^2,
$$
we have that $| m_n - \tilde{m}| > \eta $.
%\bean \label{cont} 
%\eean
As $B$ is compact, after possibly extracting a subsequence, we may assume that
$m_n $ converges to some $m^*$ in $B$, with $m^* \neq \tilde{m}$.
Since $C_n$ tends to zero, applying Theorem \ref{cv th}, there is a sequence
$m_n' $ which converges to  $\tilde{m}$
and such that
\bean  \label{prim conv}
 \int_V |{\cal C}^{-\f12}(A_{m_n'} \bh_{m_n',C_n} - \tilde{\bu})|^2  \ri 0, \quad
\bh_{m_n',C_n} \ri \tilde{\bh}
\eean
where 
$ \ds
F_{m_n', C_n} (\bh_{m_n',C_n} ) = \min_{\bg \in H^1_0(R)} 
F_{m_n', C_n} (\bg) 
$, so $F_{m_n', C_n} (\bh_{m_n',C_n} )$ converges to zero.
%, and by Theorem \ref{cv th}
%$\bh_{m_n',C_n}$ converges to $\tilde{\bh}$.
Fix $\epsilon >0$. Set 
$ \ds
F_{m_n, C_n}^{disc} (\bh_{m_n,C_n}^{disc} ) = \min_{\bg \in {\cal F}_{p_n}} 
F^{disc}_{m_n, C_n} (\bg)
$.
Let $\bh_{m_n',C_n,p}$ be the orthogonal projection of
$\bh_{m_n',C_n}$ on ${\cal F}_{p}$.
We first note that the convergence of $\bh_{m_n',C_n}$  to $\bh$
implies that $\bh_{m_n',C_n,p}$ converges to
$\bh_{m_n',C_n}$ as $p \ri \infty$, uniformly in $n$.
Thus, using minimality of $\bh_{m_n, C_n}^{disc}$,
\bea
F^{disc}_{m_n, C_n} (\bh_{m_n, C_n}^{disc})  \leq 
F^{disc}_{m_n', C_n} (\bh_{m_n',C_n,p_n}) 
%\eea
%Since $\bh_{m_n',C_n}$ is bounded and $\| \bh_{m_n',C_n,p_n}\| \leq \| \bh_{m_n',C_n}\|$,
%it follows that %$\| \bh_{m'n,C_n,p}\|$ is uniformly 
%\bea
%F^{disc}_{m_n, C_n} (\bh_{m_n',C_n,p_n}) 
\leq F^{disc}_{m_n', C_n} (\bh_{m_n',C_n}) + \epsilon,
\eea
for all $n$ large enough. Using again the boundedness of $\bh_{m_n',C_n}$, we can write  that for all
$n$ large enough,
\bea
 F^{disc}_{m_n', C_n} (\bh_{m_n',C_n}) \leq F_{m_n', C_n} (\bh_{m_n',C_n}) + \epsilon, 
\eea
and
since $F_{m_n', C_n} (\bh_{m_n',C_n} )$ converges to zero we infer that for all $n$ large enough,
 \bean \label{epsilon n}
F^{disc}_{m_n, C_n} (\bh_{m_n, C_n}^{disc})  \leq 3 \epsilon.
\eean
By Lemma \ref{bhdisc bounded}, $\bh_{m_n, C_n}^{disc}$ is bounded by a constant 
that only depends on $\tilde{\bh}$, so for all 
large enough $n$
\bean \label{epsilon 3n}
\int_V | {\cal C}^{-\f12}(A_{m^*} \bh_{m_n, C_n}^{disc} - \tilde{\bu})|^2 \leq 4 \epsilon.
\eean
 This contradicts Lemma \ref{almost stable}
for $\epsilon$ small enough.\\

\textbf{Proof of Proposition \ref{int g prop}}:\\
First we combine the exponentials in (\ref{umeas}) and (\ref{gprior})
to find
\bean 
\label{knowing m and g}
\rho(\tilde{\bu}_{meas}  | m , \bg) \rho_{{\cal F}}(\bg) \propto
\exp( -\f12 \sum_{j=1}^{N} C'(j,N) 
| ({\cal C}^{-\f12}
A_m \bg- \tilde{\bu}_{meas})(P_j) |^2 -\f12 C \int_R |\nabla\bg|^2),
\eean
which needs to be integrated in $\bg$ over ${\cal F}_p$.
With $ \bh^{disc}_{m,C}$ is as in (\ref{define bhdiscmc}) 
and the adjoint defined as in the statement of Proposition  \ref{int g prop}, 
 $\bh^{disc}_{m,C}$ satisfies
\bea
  A_m' {\cal D}^2 A_m \bh^{disc}_{m,C}+ C \bh^{disc}_{m,C}=
A_m' {\cal D}^2\tilde{\bu}_{meas}.
\eea
Setting $\bg = \bh^{disc}_{m,C}+ \bh $, it follows that
\bea
\| 
 {\cal D} A_m \bg-  {\cal D} \tilde{\bu}_{meas} \|^2 + C  \|\bg\|^2 =
\| 
 {\cal D} A_m \bh^{disc}_{m,C}-  {\cal D} \tilde{\bu}  \|^2 + C \|\bh^{disc}_{m,C}\|^2 
+ \|  {\cal D} A_m \bh \|^2 + C \|\bh\|^2,
\eea
Next we set $q= \dim {\cal F}_p$ and we introduce an orthonormal basis $\bev_1, ..., \bev_q$ 
of ${\cal F}_p$ which diagonalizes $A_m' {\cal D}^2 A_m $. Let $\mu_j^2 $ be such that 
$A_m' {\cal D}^2 A_m \bev_j = \mu_j^2 \bev_j$. We can now integrate
$\exp ( - \f12  \| 
{\cal D} A_m \bh \|^2 - \f12  C \|\bh\|^2)$ for $\bh $ over ${\cal F}_p$
by just rotating the natural basis of ${\cal F}_p$ to
the orthonormal basis $\bev_1, ..., \bev_q$ to obtain
$$
  \prod_{j=1}^q \int_{-\infty}^{\infty}
	\exp(- \f12 (\mu_j^2 + C) t^2) dt = \prod_{j=1}^q  \f{\s{2 \pi}}{\s{\mu_j^2 + C}}
	= \f{1 }{\s{\det ((2\pi)^{-1} (A_m'{\cal D}^2A_m + C I_q))}}.
 $$

\textbf{Proof of Lemma \ref{C exists}}:\\
Due to the minimization property (\ref{disc pb}) it is clear that 
$ \| {\cal D}(Ag^{(p)}  - u^{(3N)}) \| \leq \| {\cal D} u^{(3N)} \|$ for all $C>0$.
As ${\cal D}(Ag^{(p)} - v^{(3N))}$ is orthogonal to ${\cal D}(u^{(3N)} - v^{(3N)})$, by the Pythagorean theorem,
\bean
\label{pytha}
\| {\cal D}(Ag^{(p)} - u^{(3N)})  \|^2=
\| {\cal D}(Ag^{(p)} - v^{(3N)})  \|^2 + \|{\cal D}(u^{(3N)} - v^{(3N)})\|^2,
\eean
so $\| {\cal D}(Ag^{(p)}  - u^{(3N)}) \| \geq \| {\cal D}(u^{(3N)} - v^{(3N)})\| $. If we assume that
$ \|{\cal D}( Ag^{(p)}  - u^{(3N)}) \| =\| {\cal D} u^{(3N)} \|$ for some $C>0$, then 
the minimum of  (\ref{disc pb}) is achieved for $g^{(p)}=0$, so $A'{\cal D}^2u^{(3N)}=0 $
due to  (\ref{disc opt}) which contradicts 
the assumption that $v^{(3N)}$ is non-zero.
If we assume that 
$\| {\cal D}(Ag^{(p)}  - u^{(3N)}) \| =\| {\cal D}(u^{(3N)} - v^{(3N)})\| $ for some $C>0$, 
then $\|{\cal D}( Ag^{(p)}  - v^{(3N)}) \| =0 $ due to the Pythagorean theorem, 
so $C(D'D + E'E) g^{(p)} = A'{\cal D}^2(u^{(3N)} - v^{(3N)})$,
but $A'{\cal D}^2(u^{(3N)} - v^{(3N)})$ by definition of $v^{(3N)}$,
thus $g^{(p)}=0$ due to (\ref{disc opt}), leading to a contradiction. \\
In
\cite{volkov2017reconstruction}, Appendix B, we showed how $D$ and $E$ can be chosen
assuming that we use a regular grid on $R$. 
For that particular choice, $\|D\|$ and $\|E\|$ are equal to 2 while 
$\|D^{-1}\|$ and $\|E^{-1}\|$ are bounded by $\sqrt{q}$.
(we used $q$ for an upper bound for $\|D\|$ and $\|E\|$
but that bound can be improved to $\sqrt{q}$ by observing that the block matrix $M$ defined
in appendix B of \cite{volkov2017reconstruction} is the sum of the identity and a $m$- nilpotent matrix with norm 1, where $m = \sqrt{q}$).
As for any $x$ in $\RR^q$
\bean \label{coer in p}
 x^T (A'{\cal D}^2A g^{(p)} + C (D' D  + E' E ))x \geq C\|D x\|^2 + C\| E x \|^2 
\geq \f{4C}{q} \| x\|^2,
\eean
$(A'{\cal D}^2A + C (D' D  + E' E ))^{-1}$
exists for all $C>0$ and is a continuous function of $C$.
Since $g^{(p)}$ solves (\ref{disc opt}),  $g^{(p)}$ and 
$\| {\cal D}(Ag^{(p)} - v^{(3N)})\|$ are also continuous functions of $C$ in $(0, \infty)$.\\
Left multiplying  (\ref{disc opt}) by $g^p$ and applying the 
Cauchy Schwartz inequality
we find
\bea
\|{\cal D} A g^{(p)}\|^2 + C \| D g^{(p)} \|^2 +  C\|E g^{(p)}\|^2 
\leq \|{\cal D}u^{(3N)} \| \| {\cal D} A g^{(p)}\|,
\eea
thus 
\bea
\f12 \|{\cal D} A g^{(p)}\|^2 + C \| D g^{(p)} \|^2 +  C\|E g^{(p)}\|^2 
\leq \f12 \|{\cal D} u^{(3N)} \|^2 .
\eea
Recalling (\ref{coer in p}) we find
\bea
\|  g^{(p)} \|^2  \leq \f{q}{8 C} \|{\cal D} u^{(3N)} \|^2,
\eea
thus $\ds \lim_{C \ri \infty} g^{(p)}   =0$, so 
$\ds \lim_{C \ri \infty} \| {\cal D} (Ag^{(p)}  - u^{(3N)}) \| = \| {\cal D} u^{(3N)} \|$.\\
To find the limit of $Ag^{(p)} - u^{(3N)}$ as $C$ tends to zero we first recall that 
$A'{\cal D} ^2 ( u^{(3N) } - v^{(3N)})=0$.
By definition there is an $x$ in $\RR^p$ such that
${\cal D} Ax = {\cal D}  v^{(3N)}$.
From (\ref{disc opt}),
\bea
A'{\cal D} ^2A (g^{(p)} -x) + C (D' D  + E' E )(g^{(p)} -x)= -C (D' D  + E' E )x
\eea
thus
\bea
\| {\cal D} A (g^{(p)} -x) \|^2 + \f12 C \| D (g^{(p)} -x) \|^2 
+ \f12 C \| E (g^{(p)} -x) \|^2 \leq \f12 C\| D  x \|^2 
+ \f12 C \| E  x \|^2,
\eea
so $\ds \lim_{C \ri 0} \| {\cal D}A  (g^{(p)} -x)  \| = 0$ and by the Pythagorean formula
(\ref{pytha})
$$\ds \lim_{C \ri 0} \| {\cal D} (Ag^{(p)}  - u^{(3N)}) \| = \|{\cal D} ( u^{(3N)} - v^{(3N)})\|. $$

%\bibitem{S} E. P. Stephan, A boundary integral equation method for three-dimensional crack problems in elasticity, Mathematical Methods In The Applied Sciences, 1986 ,Volume: 8 Issue: 4.

\end{document}